\documentclass[10pt]{article}
\usepackage{lmodern}
\usepackage{amsmath}
\usepackage[T1]{fontenc}
\usepackage[utf8]{inputenc}
\usepackage{authblk}
\usepackage{amsfonts}
\usepackage{graphicx}
\usepackage{rotating}
\usepackage{amssymb}
\usepackage[english]{babel}
\usepackage{color}
\usepackage{amsthm}
\usepackage{graphicx}
\usepackage{mathrsfs}
\usepackage{makecell}
\usepackage{microtype}
\usepackage{mathscinet}
\usepackage{array}
\usepackage{multirow}
\usepackage{enumerate}
\usepackage[cal=boondoxo,bb=ams]{mathalfa}
\usepackage{hyperref}
\hypersetup{hidelinks}
\usepackage{booktabs}
\usepackage{arydshln}

\newtheorem{theorem}{Theorem}[section]
\newtheorem{prop}{Proposition}[section]

\newtheorem{coro}{Corollary}[section]
\newtheorem{remark}{Remark}[section]

\newcommand{\ml}{\mathcal}
\newcommand{\mb}{\mathbb}
\DeclareMathOperator{\diag}{diag}
\DeclareMathOperator{\intt}{int}
\DeclareMathOperator{\extt}{ext}
\DeclareMathOperator{\bdd}{bdd}

\def\XXint#1#2#3{{\setbox0=\hbox{$#1{#2#3}{\int}$ }
		\vcenter{\hbox{$#2#3$ }}\kern-.6\wd0}}

\title{Large time behavior for the hyperbolic-parabolic coupled system with the regularity-loss structure}
\author[1]{Wenhui Chen\thanks{Wenhui Chen (wenhui.chen.math@gmail.com)}}
\affil[1]{School of Mathematics and Information Science, Guangzhou University, \authorcr 510006 Guangzhou, P.R.  China}
\author[2]{Yan Liu\thanks{Yan Liu (ly801221@163.com)}}
\affil[2]{Department of Applied Mathematics, Guangdong University of Finance,\authorcr 510521 Guangzhou, P.R. China}
\date{}

\begin{document}
		\maketitle	
		
		\begin{abstract}
			\medskip
This paper considers the hyperbolic-parabolic coupled system, arising from the generalized thermoelastic coupled system, in the whole space $\mathbb{R}^n$. We study some qualitative properties for an energy term by diagonalization procedures, and for the solution by the WKB analysis. Particularly, we derive new large time asymptotic profiles  with the regularity-loss structure (from the biharmonic parabolic equation and the diffusion wave equation with the Riesz potential operator) and optimal decay estimates with suitable higher regularities for the Cauchy data. Finally, we discover that the wave equation with the Riesz potential dissipation is a large time approximated model of our hyperbolic-parabolic coupled system.
			\\
			
			
			\noindent\textbf{AMS Classification (2020)} 35G40, 74A15, 35B40, 	35B65
		\end{abstract}
\fontsize{12}{15}
\selectfont

\section{Introduction}\label{Section-Introduction}
\hspace{5mm}As is well-known, the reciprocal actions between elastic stresses and thermal behavior including the temperature difference of an elastic heat conductive media are modeled by thermoelastic coupled systems mathematically (cf. \cite{Green-Naghdi=1991,Chandrasekharaiah=1998,Jiang-Racke=2000}). In the last forty years, some qualitative properties (containing well-posedness, blow-up, propagation of singularity, decay property and asymptotic behavior) of solutions to some classical thermoelastic systems in bounded or unbounded domains have caught a lot of attention. We refer the interested reader to \cite{Dafermos=1968,Racke-1987,Racke-Shibata=1991,Racke-Wang=1998,Reissig-Wang=1999,Jiang-Racke=2000} and references given therein for  thermoelastic systems coupled with various laws of the heat conduction, e.g. the Fourier law, the Maxwell-Cattaneo law, or the Coleman-Gurtin law. Their main motivation is to investigate some influence from different heat conductions from parabolic equations on the well-studied (elastic) wave equations. For instance, the recent paper \cite{Chen-Takeda=2023-JDE} discovered the wave's influence as well as the heat's influence from the Fourier law on large time behavior for the classical thermoelasticity in $\mb{R}^n$. With the same philosophy as the above, our work focuses on the generalized thermoelastic system of heat equation and wave equation by the lower order coupling term.

In this manuscript, we consider the Cauchy problem for the hyperbolic-parabolic coupled system, arising from the generalized thermoelastic system (cf. \cite{Chen-Triggiani=1989,Hao-Liu=2013,Hao-Liu-Yong=2015,Fer-Liu-Racke=2019} and the explanation introduced later)  with the elastic displacement $u=u(t,x)$ and the temperature $v=v(t,x)$, as follows:
\begin{align}\label{Linear-Hyper-Para-Syst}
\begin{cases}
u_{tt}-\Delta u+v=0,&x\in\mb{R}^n,\ t>0,\\
v_t-\Delta v-u_t=0,&x\in\mb{R}^n,\ t>0,\\
(u,u_t,v)(0,x)=(u_0,u_1,v_0)(x),&x\in\mb{R}^n,
\end{cases}
\end{align}
for any $n\geqslant 1$. Our main purpose is to derive large time behavior, i.e. optimal estimates as well as asymptotic profiles, for an energy term (involving $|D|u,u_t,v$) and the solution $u$ to the Cauchy problem \eqref{Linear-Hyper-Para-Syst} by different approaches.

The coupling structures in classical thermoelastic systems based on (elastic) wave equations always imply an exponential stability, e.g. an exponential decay estimate holds for large frequencies. In fact, some physical models which posses the regularity-loss structure do not have any exponential stability, for example, the dissipative Timoshenko system \cite{Ide-Haramoto-Kawashima-2008}, the Euler-Maxwell system \cite{Hosono-Kawashima=2006}, and the Vlasov–Maxwell–Boltzmann system \cite{Duan=2011}. We will show that the regularity-loss phenomenon still holds for the thermoelastic model \eqref{Linear-Hyper-Para-Syst}, which can be regarded as the thermal effect from the heat equation $v_t-\Delta v=0$ acting on the classical wave equation $u_{tt}-\Delta u=0$ by the lower order coupling term as follows:
\begin{align*}
(+v,-u_t)^{\mathrm{T}}.
\end{align*}
 In other words, the weakly dissipative wave equation arises. The further purpose of this paper is to comprehend the weak dissipation of our model in the view of large time asymptotic profiles or approximations.

Let us recall some recent studies related to our model \eqref{Linear-Hyper-Para-Syst}. First of all,  this is a special case of the generalized thermoelastic system, or the so-called $\alpha-\beta$ coupled system, (see \cite{Chen-Triggiani=1989,Hao-Liu=2013,Hao-Liu-Yong=2015,Fer-Liu-Racke=2019} and references therein)
\begin{align}\label{alpha-beta}
\begin{cases}
u_{tt}-\Delta u-\gamma(-\Delta)^{\alpha}v=0,\\
v_t+(-\Delta)^{\beta}v+\gamma(-\Delta)^{\alpha}u_t=0,
\end{cases}
\end{align}
with $(\alpha,\beta)\in[0,1]\times[0,1]$ and $\gamma\in\mb{R}\backslash\{0\}$, namely, our model \eqref{Linear-Hyper-Para-Syst} coincides with the system \eqref{alpha-beta} when $(\alpha,\beta)=(0,1)$ and $\gamma=-1$. Note that the generalized model \eqref{alpha-beta} describes a few physical models, including the classical thermoelastic system of second order, the linear viscoelastic equation, the thermoelastic plate system (if we replace $-\Delta$ by $\Delta^2$), the Moore-Gibson-Thompson equation in thermally relaxing fluids. Moreover, the regularity analysis of \eqref{alpha-beta} is well-developed, where depending on the value of $(\alpha,\beta)$ its semigroup is classified by analytic, Gevrey class with specific order and non-smoothing. The Cauchy problem \eqref{Linear-Hyper-Para-Syst} actually was investigated by \cite{Ueda-2018}. Its author applied energy methods in the Fourier space to get the following regularity-loss type decay estimate of the energy term $\mb{E}=\mb{E}(t,x)$ such that $\mb{E}:=(\nabla u,u_t,v)^{\mathrm{T}}\in\mb{R}^{n+2}$:
\begin{align}\label{Ueda-Est}
\|\partial_x^{k}(\nabla u,u_t,v)(t,\cdot)\|_{(L^2)^{n}\times L^2\times L^2}&\lesssim (1+t)^{-\frac{n}{8}-\frac{k}{4}}\|(\nabla u_0,u_1,v_0)\|_{(L^1)^{n}\times L^1\times L^1}\notag\\
&\quad +(1+t)^{-\frac{\ell}{2}}\|\partial_x^{k+\ell}(\nabla u_0,u_1,v_0)\|_{(L^2)^{n}\times L^2\times L^2}
\end{align}
for $k,\ell\geqslant0$, and asymptotic expansions of eigenvalues to indicate the sharpness of the last estimate. Additionally, the large time approximation of $\mb{E}$ and its profile for the one dimensional model were obtained by the spectral theory with asymptotic expansions of eigenprojections. Nevertheless, optimal estimates (i.e. the sharp upper bound and lower bound estimates for large time $t\gg1$) and asymptotic profiles, particularly in the higher dimensional case, for the solution $u$, even for an energy term, are still unknown. We underline that the sharp estimate for $u$ cannot be derived from \eqref{Ueda-Est}  directly due to the fact that the energy term $\mb{E}$ does not involve the solution $u$ (see Remark \ref{Remark-Explain-Energy-u} for a more detailed explanation). For the corresponding boundary value problem to \eqref{Linear-Hyper-Para-Syst}, the polynomial strong stability and its decay rate were proposed by \cite{Ammar-B-T=1997,Liu-Rao=2005}. Another motivation for studying some qualitative properties of the solution $u$ arises from future works on 
a series of corresponding nonlinear models to the Cauchy problem \eqref{Linear-Hyper-Para-Syst} with the nonlinearity $F(u)$ on the first equation, where the authors of \cite{Clark-Lima=1997,Liu-Su=2009,Raposo-Ribeiro-Cattai=2018,Ding-Zhou=2022,Nguyen-Nhan-Truong=2023} considered global (in time) existence and blow-up of solutions in the cases $F(u)=|u|^{p-1}u$ or $F(u)=|u|^{p(x)-2}u$ under some conditions on the power exponents. We expect the detailed description of $u$ may improve some results.

In this paper, we are going to derive some large time qualitative properties for the hyperbolic-parabolic coupled system \eqref{Linear-Hyper-Para-Syst}. We first are going to adapt the diagonalization procedure developed by  \cite{Yagdjian=1997,Reissig-Wang=2005,Yang-Wang=2006,Jachmann=2008} to the model \eqref{Linear-Hyper-Para-Syst}, which allows us to derive an asymptotic representation of the energy term $W$ defined in \eqref{Energy-term}. That is the key to generalize the results in \cite{Ueda-2018} for higher dimensions. For another, by reducing the coupled system \eqref{Linear-Hyper-Para-Syst} to the higher order evolution equation, we may deduce the sharp representation of the solution $u$ in the Fourier space, in which the non-summable singularity for small frequencies can be compensated. Finally, by using the WKB analysis and the Fourier analysis, we derive some regularity-loss type decay estimates (especially, the optimal estimates under higher regularities for the Cauchy data) and large time asymptotic profiles for the energy term $W$ and the solution $u$. These profiles consist of the solutions to the biharmonic parabolic equation and the diffusion wave equation with the Riesz potential operator. As a byproduct, motivated by the singular limit process, under some conditions for the initial data, we discover the coupled system \eqref{Linear-Hyper-Para-Syst} can be approximated by the wave equation with the Riesz potential dissipation (see Appendix \ref{Appendix-A} for its large time behavior). It is a new explanation for the regularity-loss structure of our model \eqref{Linear-Hyper-Para-Syst} as well as the role of thermal equation.

\medskip
\noindent\textbf{Notation } We write $f\lesssim g$ when there  exists a positive constant $C$ such that $f\leqslant Cg$. Then, the asymptotic relation $f\simeq g$ holds if and only if $g\lesssim f\lesssim g$. The mean of a summable function $f_0=f_0(x)$ is denoted by
\begin{align*}
P_{f_0}:=\int_{\mb{R}^n}f_0(x)\mathrm{d}x.
\end{align*} We denote the identity matrix in three dimensions by $I_{3}:=\diag(1,1,1)$, the zero matrix in three dimensions by $0_{3\times 3}:=\diag(0,0,0)$, and the commutator of two matrices by $[A,B]:=AB-BA$.  Additionally, we introduce the following zones in the Fourier space:
\begin{align*}
	\ml{Z}_{\intt}(\varepsilon_0):=\{|\xi|\leqslant\varepsilon_0\ll1\},\ \ 
	\ml{Z}_{\bdd}(\varepsilon_0,N_0):=\{\varepsilon_0\leqslant |\xi|\leqslant N_0\},\ \  
	\ml{Z}_{\extt}(N_0):=\{ |\xi|\geqslant N_0\gg1\}.
\end{align*}
The smooth cut-off functions $\chi_{\intt}(\xi),\chi_{\bdd}(\xi),\chi_{\extt}(\xi)$ have the supports in their corresponding zones $\ml{Z}_{\intt}(\varepsilon_0)$, $\ml{Z}_{\bdd}(\varepsilon_0/2,2N_0)$ and $\ml{Z}_{\extt}(N_0)$, respectively, such that
\begin{align*}
	\chi_{\bdd}(\xi)=1-\chi_{\intt}(\xi)-\chi_{\extt}(\xi).
\end{align*}
The differential operator $|D|^{\sigma}$ has its symbol $|\xi|^{\sigma}$ so that $\ml{F}(|D|^{\sigma}f_0(x))=|\xi|^{\sigma}\widehat{f}_0(\xi)$ with $\sigma\in\mb{R}$, and similarly, another differential operator $\langle D\rangle^{\sigma}$ has its symbol $\langle \xi\rangle^{\sigma}$ with the Japanese bracket such that $\langle\xi\rangle^2:=1+|\xi|^2$.

\section{Main results and discussions}\setcounter{equation}{0}
\subsection{Decay property and asymptotic profile for the energy term}\hspace{5mm}Before stating our first result on the Cauchy problem \eqref{Linear-Hyper-Para-Syst}, motivated by the structure of wave operator as well as heat operator, let us introduce the energy term $W=W(t,x)\in\mb{R}^3$ such that
\begin{align}\label{Energy-term}
W:=\left(u_t+\sqrt{-\Delta}\,u,u_t-\sqrt{-\Delta}\,u,v\right)^{\mathrm{T}},
\end{align}
where $\sqrt{-\Delta}:=|D|$ owns its symbol $|\xi|$. To describe its large time asymptotic profile, we now consider two first order (in time) reference systems. Hereafter, the superscript and the subscript with the notation ``$s$'' [resp. ``$l$''] always stand for the elements derived from the small [resp. large] frequencies part.
\begin{description}
	\item[Reference System I:] Let us introduce the following biharmonic
	parabolic system equipping lower order terms (cf. \cite{Novaga-Okabe=2015,Novaga-Okabe=2016}) with its solution $V^{(s)}=V^{(s)}(t,x)\in\mb{R}^3$:
	\begin{align}\label{Reference-System-I}
		\begin{cases}
			V_t^{(s)}+\Lambda_4^{(s)}\Delta^2V^{(s)}-\Lambda_2^{(s)}\Delta V^{(s)}+\Lambda_0^{(s)}V^{(s)}=0,&x\in\mb{R}^n,\ t>0,\\
			V^{(s)}(0,x)=V^{(s)}_0(x):=Q_s^{-1}W_0(x),&x\in\mb{R}^n,
		\end{cases}
	\end{align}
	where the diagonal coefficient matrices are derived in Subsection \ref{Sub-Section-Small-Freq}, precisely,
	\begin{align*}
	\Lambda_0^{(s)}:=\diag(0,i,-i),\ \ \Lambda_2^{(s)}:=\frac{1}{2}\diag(0,1+i,1-i),\ \ \Lambda_4^{(s)}:=\frac{1}{4}\diag(4,-2-i,-2+i),
	\end{align*}
and the auxiliary matrix $Q_s=Q_s(|D|)$ with differential operators acting on the Cauchy data is defined by
\begin{align*}
Q_s:=\ml{F}^{-1}_{\xi\to x}\big(T_{1,s}(I_3+N_{2,s})(I_3+N_{3,s})(I_3+N_{4,s})\big)
\end{align*}
with the matrices $T_{1,s}$ and $N_{k,s}$ with $k=2,\dots,4$ defined in \eqref{T1-s}-\eqref{T4-s}, respectively.
	\item[Reference System II:] Let us introduce the following diffusion wave system involving the additional Riesz potential operator with its solution $V^{(l)}=V^{(l)}(t,x)\in\mb{R}^3$:
	\begin{align}\label{Reference-System-II}
	\begin{cases}
	V^{(l)}_t-\Lambda_0^{(l)}\Delta V^{(l)}+\Lambda_1^{(l)}\sqrt{-\Delta}\,V^{(l)}+\Lambda_4^{(l)}(-\Delta)^{-1}V^{(l)}=0,&x\in\mb{R}^n,\ t>0,\\
	V^{(l)}(0,x)=V_0^{(l)}(x):=Q_{l}^{-1}W_0(x),&x\in\mb{R}^n,
	\end{cases}
	\end{align}
where the diagonal coefficient matrices are derived in Subsection \ref{Sub-Section-Large-Freq}, precisely,
\begin{align*}
\Lambda_0^{(l)}:=\diag(0,0,1),\ \ \Lambda_1^{(l)}:=\diag(i,-i,0),\ \ \Lambda_4^{(l)}:=\frac{1}{2}\diag(1,1,-2),
\end{align*}
and the auxiliary matrix $Q_l=Q_l(|D|)$ with differential operators acting on the Cauchy data is defined by
\begin{align*}
	Q_l:=\ml{F}^{-1}_{\xi\to x}\big(T_{1,l}(I_3+N_{2,l})(I_3+N_{3,l})(I_3+N_{4,l})\big)
\end{align*}
with the matrices $T_{1,l}$ and $N_{k,l}$ with $k=2,\dots,4$ defined in \eqref{T1-l}-\eqref{T4-l}, respectively.
\end{description}

\begin{theorem}\label{Thm-Energy-Term}
Let $\sigma\geqslant 0$ and $\ell\geqslant0$. Suppose that the Cauchy data $W_0\in H^{\sigma+\ell}\cap L^1$ for the hyperbolic-parabolic coupled system \eqref{Linear-Hyper-Para-Syst}. Then, the energy term satisfies the following decay estimate of the regularity-loss type:
	\begin{align}\label{Upper-Bound-Est}
		\|W(t,\cdot)\|_{\dot{H}^{\sigma}}\lesssim (1+t)^{-\min\left\{\frac{n+2\sigma}{8},\frac{\ell}{2} \right\}}\|W_0\|_{H^{\sigma+\ell}\cap L^1}.
	\end{align}
Additionally, by assuming the suitable higher Sobolev regularity $\ell>\frac{n}{4}+\frac{\sigma}{2}$ and $|P_{W_0}|\neq0$, the energy term satisfies the following optimal estimate:
\begin{align}\label{Optimal-Est-W}
	t^{-\frac{n+2\sigma}{8}}|P_{W_0}|\lesssim\|W(t,\cdot)\|_{\dot{H}^{\sigma}}\lesssim t^{-\frac{n+2\sigma}{8}}\|W_0\|_{H^{\sigma+\ell}\cap L^1}
\end{align}
for large time $t\gg1$.
Furthermore, considering the lower regular  Cauchy data $W_0\in H^{\sigma+\ell-\frac{1}{2}}\cap L^1$, the energy term satisfies the following refined estimate:
\begin{align}\label{Error-ALL}
	&\left\|W(t,\cdot)-\chi_{\intt}(D)Q_sV^{(s)}(t,\cdot)-\chi_{\extt}(D)Q_lV^{(l)}(t,\cdot)\right\|_{\dot{H}^{\sigma}}\notag\\
	&\lesssim (1+t)^{-\min\left\{\frac{n+2\sigma}{8},\frac{\ell}{2} \right\}-\frac{1}{4}}\|W_0\|_{H^{\sigma+\ell-\frac{1}{2}}\cap L^1}.
\end{align}
\end{theorem}
\begin{remark}
We find the decay estimate \eqref{Upper-Bound-Est} for the energy term involving $|D|u,u_t,v$ of the regularity-loss type, namely, we have to assume an additional higher regularity $H^{\ell}$ for the Cauchy data $W_0$ to deduce our decay estimate for $W$ with the rate $(1+t)^{-\frac{\ell}{2}}$. This result exactly coincides with the estimate in \cite[Theorem 1.2]{Ueda-2018} although different energy terms are took, whose reason and comparison are referred to Remark \ref{Remark-Point-wise-without-Ueda}. Nevertheless, we do not use energy methods in the Fourier space as \cite{Ueda-2018}.
\end{remark}

\begin{remark}
As is well-known, the energy of the free wave equation $u^{\mathrm{W}}_{tt}-\Delta u^{\mathrm{W}}=0$ with the solution $u^{\mathrm{W}}=u^{\mathrm{W}}(t,x)$ ensures the conservation
\begin{align*}
\|u^{\mathrm{W}}_t(t,\cdot)\|_{L^2}+\|\,|D|u^{\mathrm{W}}(t,\cdot)\|_{L^2}=\|u^{\mathrm{W}}_t(0,\cdot)\|_{L^2}+\|\,|D|u^{\mathrm{W}}(0,\cdot)\|_{L^2},
\end{align*}
moreover, the energy of the classical heat equation $v^{\mathrm{H}}_t-\Delta v^{\mathrm{H}}=0$ with the solution $v^{\mathrm{H}}=v^{\mathrm{H}}(t,x)$ provides the polynomially decay estimate (cf. \cite[Theorem 12.1.4]{Ebert-Reissig-book}) with the Cauchy data belonging to the Bessel potential space based on the $L^1$ space, i.e. $f_0\in H^{\frac{n}{2},1}$ if and only if $\langle D\rangle^{\frac{n}{2}}f_0\in L^1$,
\begin{align*}
\|v^{\mathrm{H}}(t,\cdot)\|_{L^2}\lesssim(1+t)^{-\frac{n}{4}}\|v^{\mathrm{H}}(0,\cdot)\|_{H^{\frac{n}{2},1}}.
\end{align*}
Surprisingly, the stability of their coupled equations by the lower order term $(+v^{\mathrm{H}},-u_t^{\mathrm{W}})^{\mathrm{T}}$ is completely different from theirs. To be specific, the energy term $W(t,\cdot)$  in the $L^2$ norm to the wave-heat coupled system \eqref{Linear-Hyper-Para-Syst} decays with the polynomial rate $(1+t)^{-\frac{n}{8}}$ for the $L^1$ integrable Cauchy data, and the regularity-loss structure arises. This phenomenon never happens in classical thermoelastic systems based on (elastic) wave equations. Therefore, the way of coupling is extremely important.
\end{remark}

\begin{remark}
Comparing with the upper bound estimate \eqref{Upper-Bound-Est} and its corresponding refined (error) estimate \eqref{Error-ALL}, the decay rate and the required Sobolev regularity for the Cauchy data can be improved by $(1+t)^{-\frac{1}{4}}$ and $H^{-\frac{1}{2}}$, respectively, by subtracting the following function, i.e. its large time profile:
\begin{align*}
\chi_{\intt}(D)Q_sV^{(s)}(t,\cdot)+\chi_{\extt}(D)Q_lV^{(l)}(t,\cdot)
\end{align*}
in the $\dot{H}^{\sigma}$ norm. Consequently, we may explain the large time profile of the energy term $W$ for the hyperbolic-parabolic system \eqref{Linear-Hyper-Para-Syst} by the suitable combination of the biharmonic parabolic system \eqref{Reference-System-I} and the diffusion wave system \eqref{Reference-System-II} involving the Riesz potential operator. The Riesz potential operator $(-\Delta)^{-1}$ in the reference system \eqref{Reference-System-II} is the regularity-loss mechanism, because the solution $v=v(t,x)$ of $v_t+(-\Delta)^{-1}v=0$ fulfills
\begin{align*}
v(t,x)=\mathrm{e}^{-(-\Delta)^{-1}t}v(0,x)\ \ \mbox{and}\ \ \|v(t,\cdot)\|_{L^2}\lesssim (1+t)^{-\frac{\ell}{2}}\|v(0,\cdot)\|_{H^{\ell}}
\end{align*}
with $\ell\geqslant0$.
Our result differs from the one in \cite[Theorem 5.1]{Ueda-2018} due to the consideration of different energy terms. Additionally, thanks to the three dimensional coefficient matrix $A_0+A_1|\xi|+A_2|\xi|^2$ in the first order coupled system \eqref{Fourier-First-Order} for $\xi\in\mb{R}^n$, we can generalize the one dimensional error estimates in \cite[Theorem 5.1]{Ueda-2018} to the result for any higher dimensions $n\geqslant 1$.
\end{remark}
\begin{remark}
By the same approach as the proof of Theorem \ref{Thm-Energy-Term}, taking different regular assumptions on the Cauchy data $W_0$, we actually may derive
\begin{align*}
\mbox{LHS of \eqref{Error-ALL}}\lesssim\begin{cases}
(1+t)^{-\min\left\{\frac{n+2\sigma}{8}+\frac{1}{4},\frac{\ell}{2}\right\}}\|W_0\|_{H^{\sigma+\ell-1}\cap L^1},\\
(1+t)^{-\min\left\{\frac{n+2\sigma}{8}+\frac{1}{4},\frac{\ell}{2}+\frac{1}{2}\right\}}\|W_0\|_{H^{\sigma+\ell}\cap L^1}.
\end{cases}
\end{align*}
They imply different improvements of the decay rate or the required Sobolev regularity for the Cauchy data. Our estimate \eqref{Error-ALL} shows some improvements on the decay rate as well as on the regularity, simultaneously, for any $n\geqslant 1$ and $\sigma,\ell\geqslant0$.
\end{remark}
\begin{remark}
For the case $\ell>\frac{n}{4}+\frac{\sigma}{2}$, the estimate \eqref{Optimal-Est-W} rigorously justifies the optimality of the decay estimate
\begin{align*}
\|W(t,\cdot)\|_{\dot{H}^{\sigma}}\simeq t^{-\frac{n+2\sigma}{8}}\ \ \mbox{provided}\ \ |P_{W_0}|\neq0. 
\end{align*} The question on a rigorous justification for the large time optimal estimate with some lower regular data when $0\leqslant \ell\leqslant \frac{n}{4}+\frac{\sigma}{2}$ is still open.
\end{remark}

\subsection{Decay property and asymptotic profile for the solution}\hspace{5mm}
Because the energy term $W$ defined in \eqref{Energy-term} does not contain the solution $u$, our second result turns to some qualitative properties for it, which provides a new perspective for the hyperbolic-parabolic coupled system \eqref{Linear-Hyper-Para-Syst}. We next state a note for the difficulty by using energy estimates to derive estimates for the solution $u$, which arises from the second order hyperbolic equation.
\begin{remark}\label{Remark-Explain-Energy-u}
	Let us explain that estimating the solution $u$ is not a simple generalization of an application for the energy term $\mb{E}$ in \cite{Ueda-2018}. If one uses the sharp pointwise estimate of the energy term $\mb{E}$ in the Fourier space, to be specific, \cite[Inequality (1.15)]{Ueda-2018}, we may get
	\begin{align*}
		\chi_{\intt}(\xi)|\xi|\,|\widehat{u}|\lesssim\chi_{\intt}(\xi)\mathrm{e}^{-c|\xi|^4t}(|\xi|\,|\widehat{u}_0|+|\widehat{u}_1|+|\widehat{v}_0|).
	\end{align*}
By directly dividing $|\xi|$ on the both sides of it, the next unbounded estimate of the solution for small frequencies in lower dimensions arises:
\begin{align*}
\|\chi_{\intt}(D)u(t,\cdot)\|_{L^2}&\lesssim \left\|\chi_{\intt}(\xi)\mathrm{e}^{-c|\xi|^4t}\right\|_{L^2}\|u_0\|_{L^1}+\left\|\chi_{\intt}(\xi)|\xi|^{-1}\mathrm{e}^{-c|\xi|^4t}\right\|_{L^2}\|(u_1,v_0)\|_{L^1\times L^1}\\
&\lesssim (1+t)^{-\frac{n}{8}}\|u_0\|_{L^1}+\left(\int_0^{\varepsilon_0}r^{n-3}\mathrm{e}^{-2cr^4t}\mathrm{d}r\right)^{1/2}\|(u_1,v_0)\|_{L^1\times L^1},
\end{align*}
because of the non-summable singularity at $r=0$ when $n=1,2$. This phenomenon is caused by some information hid in the total energy term $\mb{E}$. The singularity $\chi_{\intt}(\xi)|\xi|^{-1}$ can be compensated if we study $\widehat{u}$ independently in a deep way by another manner.
\end{remark}
\noindent Before showing this result, let us introduce the differential operators in the semigroup form to describe an asymptotic profile
\begin{align*}
\ml{G}^{(s)}(t,|D|):=\chi_{\intt}(D)\mathrm{e}^{-\Delta^2t}\ \ \mbox{and}\ \ \ml{G}^{(l)}(t,|D|):=\chi_{\extt}(D)\cos(|D|t)\mathrm{e}^{-\frac{1}{2}(-\Delta)^{-1}t}.
\end{align*}
Note that $\ml{G}^{(s)}(t,|D|)$ is the biharmonic heat semigroup, and $\ml{G}^{(l)}(t,|D|)$ consists of the half-wave semigroup $\frac{1}{2}(\mathrm{e}^{i|D|t}+\mathrm{e}^{-i|D|t})=\cos(|D|t)$ as well as the Riesz potential semigroup. 

\begin{theorem}\label{Thm-Solution-Itself}
Let $\ell\geqslant0$. Suppose that the Cauchy data $u_0\in H^{\ell}\cap L^1$, $u_1\in H^{\ell-1}\cap L^1$ and $v_0\in H^{\ell-3}\cap L^1$  for the hyperbolic-parabolic coupled system \eqref{Linear-Hyper-Para-Syst}. Then, the solution satisfies the following decay estimate of the regularity-loss type:
\begin{align}\label{Est-05}
	\|u(t,\cdot)\|_{L^2}\lesssim (1+t)^{-\frac{n}{8}}\|(u_0,u_1,v_0)\|_{L^1\times L^1\times L^1}+(1+t)^{-\frac{\ell}{2}}\|(u_0,u_1,v_0)\|_{H^{\ell}\times H^{\ell-1}\times H^{\ell-3}}.
\end{align}
Additionally, by assuming the suitable higher Sobolev regularity $\ell>\frac{n}{4}$ and $|P_{u_0-v_0}|\neq0$, the solution satisfies the following optimal estimate:
\begin{align}\label{Est-06}
t^{-\frac{n}{8}}|P_{u_0-v_0}|\lesssim\|u(t,\cdot)\|_{L^2}\lesssim t^{-\frac{n}{8}}\|(u_0,u_1,v_0)\|_{(H^{\ell}\cap L^1)\times (H^{\ell-1}\cap L^1)\times (H^{\ell-3}\cap L^1)}
\end{align}
for large time $t\gg1$. Furthermore, considering the lower regular  Cauchy data $u_0\in H^{\ell-1}\cap L^1$, the solution satisfies the following refined estimate:
\begin{align}\label{Est-08}
	&\left\|u(t,\cdot)-\ml{G}^{(s)}(t,|D|)\big(u_0(\cdot)-v_0(\cdot)\big)-\ml{G}^{(l)}(t,|D|)u_0(\cdot)\right\|_{L^2}\notag\\
	&\lesssim (1+t)^{-\min\left\{\frac{n}{8}+\frac{1}{2},\frac{n}{4}\right\}}\|(u_0,u_1,v_0)\|_{L^1\times L^1\times L^1}+(1+t)^{-\frac{\ell}{2}}\|(u_0,u_1,v_0)\|_{H^{\ell-1}\times H^{\ell-1}\times H^{\ell-3}}.
\end{align}
\end{theorem}
\begin{remark}
The estimate \eqref{Est-05} shows the regularity-loss phenomenon occurs when $\ell>1$ due to the natural assumptions $u_0\in H^1$ and $u_1\in L^2$ for the free wave equation. In other words, the new regularity-loss threshold of the large time optimal estimate \eqref{Est-06} for the hyperbolic-parabolic coupled system \eqref{Linear-Hyper-Para-Syst} is the dimension $n=4$. When $n=1,2,3$, we just need to set 
\begin{align*}
(H^1\cap L^1)\times (L^2\cap L^1)\times (L^2\cap L^1)
\end{align*}
 for the Cauchy data.
\end{remark}
\begin{remark}\label{Remark-01}
Comparing with the upper bound estimate \eqref{Est-05} and its corresponding refined estimate \eqref{Est-08}, the decay rate for the $L^1$ Cauchy data can be improved by $(1+t)^{-\min\left\{\frac{1}{2},\frac{n}{8}\right\}}$, and the required Sobolev regularity for the initial data $u_0$ can be improved by $H^{-1}$, by subtracting the following function, i.e. its large time profile:
\begin{align*}
\chi_{\intt}(D)\mathrm{e}^{-\Delta^2t}\big(u_0(\cdot)-v_0(\cdot)\big)+\chi_{\extt}(D)\cos(|D|t)\mathrm{e}^{-\frac{1}{2}(-\Delta)^{-1}t}u_0(\cdot)
\end{align*}
in the $L^2$ norm. Consequently, we may explain the large time profile of the solution $u$ for the hyperbolic-parabolic system \eqref{Linear-Hyper-Para-Syst} by the suitable combination of biharmonic heat equation, half-wave equation involving Riesz potential operator. This is a new insight of the model.
\end{remark}

To end this section, we explain the regularity-loss type decay property for the coupled system \eqref{Linear-Hyper-Para-Syst} by the following wave equation with the Riesz potential dissipation (we call it as the weakly dissipative wave equation later):
\begin{align}\label{Wave-Riesz-Dissipation}
\begin{cases}
w_{tt}-\Delta w+(-\Delta)^{-1}w_t=0,&x\in\mb{R}^n,\ t>0,\\
w(0,x)=w_0(x),\ w_t(0,x)=w_1(x),&x\in\mb{R}^n,
\end{cases}
\end{align}
under some conditions on the initial data.  Remark that we easily derive large time behavior of the solution $w$ to the Cauchy problem \eqref{Wave-Riesz-Dissipation} in Appendix \ref{Appendix-A} by the WKB analysis and the Fourier analysis. This equation of $w=w(t,x)$ is strongly motivated by the formal singular limit of our model \eqref{Linear-Hyper-Para-Syst}, precisely, the consideration of the unknown functions $(u^{\epsilon},v^{\epsilon})=(u^{\epsilon}(t,x),v^{\epsilon}(t,x))$ with the parameter $0\leqslant \epsilon<1$ (note that $\epsilon$ is the superscript of the function instead of the power) solving
\begin{align}\label{Singular-Limits-Hyper-Para}
	\begin{cases}
 u^{\epsilon}_{tt}-\Delta u^{\epsilon}+v^{\epsilon}=0,&x\in\mb{R}^n,\ t>0,\\
		\epsilon v^{\epsilon}_t-\Delta v^{\epsilon}-u^{\epsilon}_t=0,&x\in\mb{R}^n,\ t>0,\\
		(u^{\epsilon},u^{\epsilon}_t,v^{\epsilon})(0,x)=(u^{\epsilon}_0,u^{\epsilon}_1,v^{\epsilon}_0)(x),&x\in\mb{R}^n.
	\end{cases}
\end{align}
Formally taking the singular limit $\epsilon=0$, the thermal equation \eqref{Singular-Limits-Hyper-Para}$_2$ turns into $v^0=(-\Delta)^{-1}u^0_t$. Then, the formal singular limit model of \eqref{Singular-Limits-Hyper-Para} is given by \eqref{Wave-Riesz-Dissipation}. It hints some relations between these two models, even for large time.

Summarizing the results in Theorem \ref{Thm-Solution-Itself} and Proposition \ref{Prop-Riesz-wave}, we conclude the next table for the comparison. It shows the same decay rate for the $L^1$ data, the same required regularity for the Sobolev data with the decay rate $(1+t)^{-\frac{\ell}{2}}$, and the same structure of large time profiles. But improved decay rates and improved required regularities for the first data in the error estimates between these models are slightly different due to the coupled effect.
\begin{table}[h!]
	\begin{center}
		\caption{Large time behavior for two models with the regularity-loss structure}
		\medskip
		\label{tab:table1}
		\begin{tabular}{cccc} 
			\toprule
			& Decay rate of the  & Regularity of the  & Large time  \\
			&$L^1$ data &Sobolev data &asymptotic profile\\
			\midrule
			The weakly dissipative&\multirow{2}{*}{$(1+t)^{-\frac{n}{8}}$} & \multirow{2}{*}{$H^{\ell}\times H^{\ell-1}$} & \multirow{2}{*}{$\ml{G}^{(s)}(t,|D|)w_0+\ml{G}^{(l)}(t,|D|)w_0$} \\
			wave equation \eqref{Wave-Riesz-Dissipation} & &  &\\
			\midrule
			The hyperbolic-parabolic&\multirow{2}{*}{$(1+t)^{-\frac{n}{8}}$} & \multirow{2}{*}{$H^{\ell}\times H^{\ell-1}\times H^{\ell-3}$} & \multirow{2}{*}{$\ml{G}^{(s)}(t,|D|)(u_0-v_0)+\ml{G}^{(l)}(t,|D|)u_0$} \\
			coupled system \eqref{Linear-Hyper-Para-Syst} & &  &\\
			\bottomrule
			\multicolumn{4}{l}{\emph{$*$ All decay rates in the above for the $L^2$ norm are optimal under $|P_{w_0}|\neq0\neq|P_{u_0-v_0}|$ when $\ell>\frac{n}{4}$.}}
		\end{tabular}
	\end{center}
\end{table}

Therefore, by taking $w_0=u_0$ in the weakly dissipative wave equation \eqref{Wave-Riesz-Dissipation} and $v_0=0$ in the hyperbolic-parabolic coupled system \eqref{Linear-Hyper-Para-Syst}, regarding the profile
\begin{align*}
\ml{G}^{(s)}(t,|D|)u_0+\ml{G}^{(l)}(t,|D|)u_0
\end{align*}
as the bridge, we are able to prove the next large time approximation by the triangle inequality between two error estimates from Theorem \ref{Thm-Solution-Itself} and Proposition \ref{Prop-Riesz-wave}.
\begin{coro}\label{Coro-01}
Let $\ell\geqslant0$. Suppose that the Cauchy data $u_0,u_1\in H^{\ell-1}\cap L^1$ and $v_0=0$  for the hyperbolic-parabolic coupled system \eqref{Linear-Hyper-Para-Syst}, the Cauchy data $w_0=u_0$ and $w_1\in H^{\ell-1}\cap L^1$ for the weakly dissipative wave equation \eqref{Wave-Riesz-Dissipation}. Then, their solutions satisfy the following approximation:
\begin{align*}
	\|u(t,\cdot)-w(t,\cdot)\|_{L^2}&\lesssim (1+t)^{-\min\left\{\frac{n}{8}+\frac{1}{2},\frac{n}{4}\right\}}\|(u_0,u_1,w_1)\|_{L^1\times L^1\times L^1}\notag\\
	&\quad+(1+t)^{-\frac{\ell}{2}}\|(u_0,u_1,w_1)\|_{H^{\ell-1}\times H^{\ell-1}\times H^{\ell-1}}.
\end{align*}
Particularly, when $\ell>\frac{n}{4}$, it holds
\begin{align*}
\lim\limits_{t\to+\infty}t^{\frac{n}{8}}\|u(t,\cdot)-w(t,\cdot)\|_{L^2}=0.
\end{align*}
\end{coro} 
\begin{remark} Let us recall Theorem \ref{Thm-Solution-Itself} and Proposition \ref{Prop-Riesz-wave}.
With the same reason as Remark \ref{Remark-01}, i.e. the improved decay rate for the $L^1$ data and the improved required Sobolev regularity for $u_0$, we find that the weakly dissipative wave equation \eqref{Wave-Riesz-Dissipation} is an approximated model of the hyperbolic-parabolic coupled system \eqref{Linear-Hyper-Para-Syst} for large time $t\gg1$. It also explains the phenomenon in Table \ref{tab:table1}. This is a new insight of the coupled system \eqref{Linear-Hyper-Para-Syst}
\end{remark}
\begin{remark}
Motivated by the formal singular limit in \eqref{Singular-Limits-Hyper-Para}, we conjecture that an approximated model of the generalized thermoelastic coupled system \eqref{alpha-beta} is the wave equation with the general structural damping
\begin{align*}
w_{tt}-\Delta w+\gamma^2(-\Delta)^{2\alpha-\beta}w_t=0
\end{align*}
with $(\alpha,\beta)\in[0,1]\times[0,1]$ and $\gamma\in\mb{R}\backslash\{0\}$. Our result in Corollary \ref{Coro-01} verifies this conjecture when $\alpha=0$, $\beta=1$ and $\gamma=1$. Nevertheless, the rigorous justification for general cases is beyond the scope of this manuscript.
\end{remark}

\section{Large time asymptotic behavior for the energy term}\setcounter{equation}{0}
\hspace{5mm}From the hyperbolic-parabolic coupled system \eqref{Linear-Hyper-Para-Syst}, the energy term $W$ defined in \eqref{Energy-term} satisfies the first order (in time) differential system
\begin{align*}
W_t=\begin{pmatrix}
	u_{tt}+\sqrt{-\Delta}\,u_t\\
	u_{tt}-\sqrt{-\Delta}\,u_t\\
	v_t
\end{pmatrix}=-\frac{1}{2} 
\begin{pmatrix}
	0 & -2\sqrt{-\Delta}&2\\
	2\sqrt{-\Delta} & 0 & 2\\
	-1 & -1 & -2\Delta 
\end{pmatrix}W,
\end{align*}
namely,
\begin{align*}
\begin{cases}
W_t+A_0W+A_1\sqrt{-\Delta}\, W-A_2\Delta W=0,&x\in\mb{R}^n,\ t>0,\\
W(0,x)=W_0(x),&x\in\mb{R}^n,
\end{cases}
\end{align*} 
whose coefficient matrices are given by
\begin{align*}
A_0=\frac{1}{2}\begin{pmatrix}
0 & 0 & 2\\
0 & 0 & 2\\
-1 & -1 & 0
\end{pmatrix},\quad A_1=\begin{pmatrix}
0 & -1 & 0\\
1 & 0 & 0\\
0 & 0 & 0
\end{pmatrix},\quad A_2=\begin{pmatrix}
0 & 0 & 0\\
0 & 0 & 0\\
0 & 0 & 1
\end{pmatrix}.
\end{align*}
\begin{remark}
Due to the different constructions of energy term between ours and the one in \cite[Page 682]{Ueda-2018}, the corresponding first order (in time) systems of the hyperbolic-parabolic coupled system \eqref{Linear-Hyper-Para-Syst} are quite different, especially, the coefficient matrix of the zero order term. Specifically, the matrix $A_0$ is not even anti-symmetric, which violates the condition in classical results of hyperbolic-parabolic coupled systems.
\end{remark}

Let us apply the partial Fourier transform with respect to the spatial variable, i.e. $\widehat{W}=\widehat{W}(t,\xi)$ such that $\widehat{W}=\ml{F}_{x\to\xi}(W)$ which solves the $|\xi|$-dependent first order coupled system
\begin{align}\label{Fourier-First-Order}
\begin{cases}
\widehat{W}_t+\left(A_0+A_1|\xi|+A_2|\xi|^2\right)\widehat{W}=0,&\xi\in\mb{R}^n,\ t>0,\\
\widehat{W}(0,\xi)=\widehat{W}_0(\xi),&\xi\in\mb{R}^n,
\end{cases}
\end{align}
with the Cauchy data $\widehat{W}_0=\widehat{W}_0(\xi)$ such that
\begin{align*}
\widehat{W}_0=(\widehat{u}_1+|\xi|\widehat{u}_0,\widehat{u}_1-|\xi|\widehat{u}_0,\widehat{v}_0)^{\mathrm{T}}.
\end{align*}
In order to understand the dominant coefficient influenced by the magnitude of $|\xi|$, we may divide our discussion into three cases: small frequencies $\xi\in\ml{Z}_{\intt}(\varepsilon_0)$, large frequencies $\xi\in\ml{Z}_{\extt}(N_0)$ and bounded frequencies $\xi\in\ml{Z}_{\bdd}(\varepsilon_0,N_0)$ with $0<\varepsilon_0\ll 1$ as well as $N_0\gg1$.

\subsection{Diagonalization procedure for small frequencies}\label{Sub-Section-Small-Freq}
\hspace{5mm}We first diagonalize the dominant coefficient matrix $A_0$ by taking
\begin{align}\label{T1-s}
T_{1,s}:=\begin{pmatrix}
-1 & -i & i\\
1 & -i & i\\
0 & 1 & 1
\end{pmatrix},
\end{align}
so that
\begin{align*}
	T_{1,s}^{-1}A_0T_{1,s}=:\Lambda_0^{(s)}=\diag(0,i,-i).
\end{align*}
By multiplying $T_{1,s}^{-1}$ to the equation of \eqref{Fourier-First-Order}, the new unknown $\widehat{W}^{(1,s)}:=T_{1,s}^{-1}\widehat{W}$ satisfies
\begin{align*}
\widehat{W}^{(1,s)}_t+\left(\Lambda_0^{(s)}+A_1^{(1,s)}|\xi|+A_2^{(1,s)}|\xi|^2\right)\widehat{W}^{(1,s)}=0,
\end{align*}
where its coefficient matrices are
\begin{align*}
A_1^{(1,s)}=\frac{1}{2}\begin{pmatrix}
	0 & -2i & 2i\\
	-i & 0 & 0\\
	i & 0 & 0
\end{pmatrix}\ \ \mbox{and}\ \ A_2^{(1,s)}=\frac{1}{2}\begin{pmatrix}
0 & 0 & 0\\
0 & 1 & 1\\
0 & 1 & 1
\end{pmatrix}.
\end{align*}

To retain the first diagonal matrix $\Lambda_0^{(s)}$ and deduce the next diagonal matrix from $A_1^{(1,s)}$, we now introduce
\begin{align}\label{T2-s}
T_{2,s}:=I_3+N_{2,s}\ \ \mbox{with}\ \ N_{2,s}:=\frac{|\xi|}{2}\begin{pmatrix}
0 & -2 & -2\\
1 & 0 & 0\\
1 & 0 & 0
\end{pmatrix}.
\end{align}
For this reason, the unknown $\widehat{W}^{(2,s)}:=T_{2,s}^{-1}\widehat{W}^{(1,s)}$ solves
\begin{align*}
\widehat{W}^{(2,s)}_t+T_{2,s}^{-1}\Lambda_0^{(s)}T_{2,s}\widehat{W}^{(2,s)}+T_{2,s}^{-1}\left(A_1^{(1,s)}|\xi|+A_2^{(1,s)}|\xi|^2\right)T_{2,s}\widehat{W}^{(2,s)}=0.
\end{align*}
We may separate the diagonal part in the new matrix by
\begin{align*}
T_{2,s}^{-1}\Lambda_0^{(s)}T_{2,s}=\Lambda_0^{(s)}+T_{2,s}^{-1}[\Lambda_0^{(s)},T_{2,s}]=\Lambda_0^{(s)}-T_{2,s}^{-1}[N_{2,s},\Lambda_0^{(s)}].
\end{align*}
Because of the decomposition
\begin{align*}
T_{2,s}^{-1}A_1^{(1,s)}T_{2,s}=T_{2,s}^{-1}A_1^{(1,s)}+T_{2,s}^{-1}A_1^{(1,s)}N_{2,s},
\end{align*}
we write
\begin{align*}
\widehat{W}^{(2,s)}_t&+\Lambda_0^{(s)}\widehat{W}^{(2,s)}+T_{2,s}^{-1}\left[\left(A_1^{(1,s)}|\xi|-[N_{2,s},\Lambda_0^{(s)}]\right)+\left(A_1^{(1,s)}N_{2,s}|\xi|+A_2^{(1,s)}T_{2,s}|\xi|^2\right)\right]\widehat{W}^{(2,s)}=0.
\end{align*}
The auxiliary matrix $N_{2,s}$ follows
\begin{align*}
	A_1^{(1,s)}|\xi|-[N_{2,s},\Lambda_0^{(s)}]=:\Lambda_1^{(s)}|\xi|=0_{3\times 3}\ \ \mbox{and}\ \ T_{2,s}^{-1}+N_{2,s}T_{2,s}^{-1}=I_3.
\end{align*}
Then, the first order system arises
\begin{align*}
\widehat{W}^{(2,s)}_t+\Lambda_0^{(s)}\widehat{W}^{(2,s)}+\left(A_1^{(2,s)}|\xi|^2+R_{2,s}\right)\widehat{W}^{(2,s)}=0,
\end{align*}
where to simplify our notations we denote
\begin{align*}
A_1^{(2,s)}|\xi|^2:=A_1^{(1,s)}N_{2,s}|\xi|+A_2^{(1,s)}|\xi|^2=\frac{1}{2}\begin{pmatrix}
0 & 0 & 0\\
0 & 1+i & 1+i\\
0 & 1-i & 1-i 
\end{pmatrix}|\xi|^2,
\end{align*}
and the remainder $R_{2,s}=O(|\xi|^3)$ being
\begin{align*}
R_{2,s}:=A_2^{(1,s)}N_{2,s}|\xi|^2-N_{2,s}T_{2,s}^{-1}\left(A_1^{(1,s)}N_{2,s}|\xi|+A_2^{(1,s)}T_{2,s}|\xi|^2\right).
\end{align*}

With the same reason as the last step, we introduce
\begin{align}\label{T3-s}
T_{3,s}:=I_3+N_{3,s}\ \ \mbox{with}\ \ N_{3,s}:=\frac{|\xi|^2}{4}\begin{pmatrix}
0 & 0 & 0\\
0 & 0 & -1+i\\
0 & -1-i & 0
\end{pmatrix}.
\end{align}
Thus, the unknown $\widehat{W}^{(3,s)}:=T_{3,s}^{-1}\widehat{W}^{(2,s)}$ satisfies
\begin{align*}
\widehat{W}^{(3,s)}_t+\Lambda_0^{(s)}\widehat{W}^{(3,s)}+T_{3,s}^{-1}\left[\left(A_1^{(2,s)}|\xi|^2-[N_{3,s},\Lambda_0^{(s)}]\right)+\left(A_1^{(2,s)}N_{3,s}|\xi|^2+R_{2,s}T_{3,s}\right)\right]\widehat{W}^{(3,s)}=0.
\end{align*}
Thanks to the auxiliary matrix $N_{3,s}$, it gives
\begin{align*}
A_1^{(2,s)}|\xi|^2-[N_{3,s},\Lambda_0^{(s)}]=:\Lambda_2^{(s)}|\xi|^2=\frac{1}{2}\diag(0,1+i,1-i)|\xi|^2,
\end{align*}
which leads to the following first order system:
\begin{align*}
\widehat{W}^{(3,s)}_t+\left(\Lambda_0^{(s)}+\Lambda_2^{(s)}|\xi|^2\right)\widehat{W}^{(3,s)}+R_{3,s}\widehat{W}^{(3,s)}=0.
\end{align*}
Here, the remainder $R_{3,s}=O(|\xi|^3)$ can be written by
\begin{align*}
R_{3,s}:\!&=-N_{3,s}T_{3,s}^{-1}\Lambda_2^{(s)}|\xi|^2+T_{3,s}^{-1}\left(A_1^{(2,s)}N_{3,s}|\xi|^2+R_{2,s}T_{3,s}\right)\\
& =T_{3,s}^{-1}\left(A_1^{(2,s)}|\xi|^2+R_{2,s}\right)N_{3,s}-N_{3,s}T_{3,s}^{-1}\left(\Lambda_2^{(s)}|\xi|^2+R_{2,s}\right)+R_{2,s},
\end{align*}
whose dominant term $R_{2,s}$ may be decomposed into
\begin{align*}
R_{2,s}&=[A_2^{(1,s)},N_{2,s}]|\xi|^2-N_{2,s}A_1^{(1,s)}N_{2,s}|\xi|+N_{2,s}^2T_{2,s}^{-1}\left(A_1^{(1,s)}N_{2,s}|\xi|+A_2^{(1,s)}|\xi|^2\right)\\
&\quad-N_{2,s}T_{2,s}^{-1}A_2^{(1,s)}N_{2,s}|\xi|^2.
\end{align*}

We until now have derived a first order coupled system with respect to $\widehat{W}^{(3,s)}$ containing the dominant diagonal matrix $\Lambda_0^{(s)}+\Lambda_2^{(s)}|\xi|^2$. Nevertheless, the first characteristic root is zero, which implies that our diagonalization procedure does not finish. For this reason, it requires further steps of diagonalization procedure to obtain pairwise distinct characteristic roots with non-zero real parts. Before doing this step, we should take the dominant terms in $R_{3,s}$, namely, the lowest order terms with respect to $|\xi|$ in $R_{2,s}$, that is
\begin{align*}
A_1^{(3,s)}|\xi|^3:=[A_2^{(1,s)},N_{2,s}]|\xi|^2-N_{2,s}A_1^{(1,s)}N_{2,s}|\xi|=\frac{1}{2}\begin{pmatrix}
0 & 2 & 2\\
1 & 0 & 0\\
1 & 0 & 0
\end{pmatrix}|\xi|^3.
\end{align*}
Facing with the new system
\begin{align*}
\widehat{W}^{(3,s)}_t+\left(\Lambda_0^{(s)}+\Lambda_2^{(s)}|\xi|^2\right)\widehat{W}^{(3,s)}+\left(A_1^{(3,s)}|\xi|^3+\widetilde{R}_{3,s}\right)\widehat{W}^{(3,s)}=0
\end{align*}
with the remainder $\widetilde{R}_{3,s}:=R_{3,s}-A_1^{(3,s)}|\xi|^3=O(|\xi|^4)$, we define the new unknown $\widehat{W}^{(4,s)}:=T_{4,s}^{-1}\widehat{W}^{(3,s)}$ and the auxiliary matrix
\begin{align}\label{T4-s}
T_{4,s}:=I_3+N_{4,s}\ \ \mbox{with}\ \ N_{4,s}:=\frac{|\xi|^3}{2}\begin{pmatrix}
	0 & -2i & 2i\\
	i & 0 & 0\\
	-i & 0 &0
\end{pmatrix}.
\end{align}
Due to the construction that
\begin{align*}
A_1^{(3,s)}|\xi|^3-[N_{4,s},\Lambda_0^{(s)}]=:\Lambda_3^{(s)}|\xi|^3=0_{3\times 3},
\end{align*}
similarly to the previous computations, we are able to derive
\begin{align*}
\widehat{W}^{(4,s)}_t+\left(\Lambda_0^{(s)}+\Lambda_2^{(s)}|\xi|^2\right)\widehat{W}^{(4,s)}+\left(A_1^{(4,s)}|\xi|^4+R_{4,s}\right)\widehat{W}^{(4,s)}=0,
\end{align*}
where the remainder $R_4=O(|\xi|^5)$ is
\begin{align*}
R_{4,s}&:=-N_{4,s}T_{4,s}^{-1}\Lambda_2^{(s)}|\xi|^2+T_{4,s}^{-1}\left(\Lambda_2^{(s)}|\xi|^2+A_1^{(3,s)}|\xi|^3\right)+T_{4,s}^{-1}\widetilde{R}_{3,s}N_{4,s}\\
&\ \quad -N_{4,s}T_{4,s}^{-1}\widetilde{R}_{3,s}+\left(\widetilde{R}_{3,s}-A_1^{(4,s)}|\xi|^4\right),
\end{align*}
moreover, the coefficient matrix $A_1^{(4,s)}$ is
\begin{align*}
A_1^{(4,s)}|\xi|^4:\!&=A_1^{(2,s)}N_{3,s}|\xi|^2-N_{3,s}\Lambda_2^{(s)}|\xi|^2+N_{2,s}^2A_1^{(1,s)}N_{2,s}|\xi|+N_{2,s}[N_{2,s},A_2^{(1,s)}]|\xi|^2\\
& =\frac{1}{4}\begin{pmatrix}
	4 & 0 & 0\\
	0 & -2-i & -3-i\\
	0 & -3+i & -2+i
\end{pmatrix}|\xi|^4.
\end{align*}

Finally, let us take the unknown $\widehat{W}^{(5,s)}:=T_{5,s}^{-1}\widehat{W}^{(4,s)}$, where
\begin{align}\label{T5-s}
T_{5,s}:=I_3+N_{5,s}\ \ \mbox{with}\ \ N_{5,s}:=\frac{|\xi|^4}{8}\begin{pmatrix}
	0 & 0 & 0\\
	0 & 0 & 1-3i\\
	0 & 1+3i & 0
\end{pmatrix}.
\end{align}
Due to the fact that
\begin{align*}
A_1^{(4,s)}|\xi|^4-[N_{5,s},\Lambda_0^{(s)}]=:\Lambda_4^{(s)}|\xi|^4=\frac{1}{4}\diag(4,-2-i,-2+i)|\xi|^4,
\end{align*}
we immediately derive
\begin{align*}
\widehat{W}^{(5,s)}_t+\left(\Lambda_0^{(s)}+\Lambda_2^{(s)}|\xi|^2+\Lambda_4^{(s)}|\xi|^4\right)\widehat{W}^{(5,s)}+R_{5,s}\widehat{W}^{(5,s)}=0
\end{align*}
with the suitable remainder $R_{5,s}=O(|\xi|^5)$.

To sum up, by carrying out five steps of diagonalization procedure, we have derived the diagonal matrix in the final coupled system
\begin{align*}
\Lambda_0^{(s)}+\Lambda_2^{(s)}|\xi|^2+\Lambda_4^{(s)}|\xi|^4=\diag\left(|\xi|^4,i+\frac{1+i}{2}|\xi|^2-\frac{2+i}{4}|\xi|^4,-i+\frac{1-i}{2}|\xi|^2-\frac{2-i}{4}|\xi|^4\right).
\end{align*}
Thanks to the positive real parts as well as the pairwise distinct diagonal elements of the last matrix and the higher order term $R_{5,s}$ to be the remainder, we in Proposition \ref{Prop-Diag-Small} successfully obtain an asymptotic representation of the energy term in the Fourier space for small frequencies according to the general philosophy proposed in \cite[Chapter 2]{Jachmann=2008}.
\begin{prop}\label{Prop-Diag-Small}
The characteristic roots $\lambda_{1,2,3}$ of the coefficient matrix $A_0+A_1|\xi|+A_2|\xi|^2$ in the Cauchy problem \eqref{Fourier-First-Order} behave for $\xi\in\ml{Z}_{\intt}(\varepsilon_0)$ as
\begin{align*}
\lambda_1=-|\xi|^4+O(|\xi|^5),\ \ \lambda_{2,3}=\mp i-\frac{1\pm i}{2}|\xi|^2+\frac{2\pm i}{4}|\xi|^4+O(|\xi|^5).
\end{align*}
Moreover, the solution for $\xi\in\ml{Z}_{\intt}(\varepsilon_0)$  to the Cauchy problem \eqref{Fourier-First-Order} has the following representation: 
\begin{align*}
\chi_{\intt}(\xi)\widehat{W}=\chi_{\intt}(\xi)T_s\diag\left(\mathrm{e}^{\lambda_1 t},\mathrm{e}^{\lambda_2 t},\mathrm{e}^{\lambda_3 t}\right)T_s^{-1}\widehat{W}_0
\end{align*}
with the coefficient $T_s=O(1)$ such that
\begin{align*}
T_s:=T_{1,s}(I_3+N_{2,s})(I_3+N_{3,s})(I_3+N_{4,s})(I_3+N_{5,s}),
\end{align*}
where the matrices $T_{1,s}$ and $N_{k,s}$ with $k=2,\dots,5$ are defined in \eqref{T1-s}-\eqref{T5-s}, respectively.
\end{prop}

\subsection{Diagonalization procedure for large frequencies}\label{Sub-Section-Large-Freq}
\hspace{5mm}Notice that the dominant coefficient matrix $A_2$ is diagonal already, namely,
\begin{align*}
\Lambda_0^{(l)}:=A_2=\diag(0,0,1).
\end{align*}

Because of the identity of the first two diagonal elements of $\Lambda_0^{(l)}$, our standard diagonalization procedure used in the last subsection does not work now. To overcome the technical difficulty, let us introduce the unknown $\widehat{W}^{(1,l)}:=T_{1,l}^{-1}\widehat{W}$ with the block diagonal matrix
\begin{align}\label{T1-l}
T_{1,l}:=\left(\begin{array}{cc;{2pt/2pt}c}
	i & -i & 0\\
	1 & 1 & 0\\ \hdashline[2pt/2pt] 
	0 & 0 & 1
\end{array}\right)
\end{align}
so that
\begin{align*}
T_{1,l}^{-1}A_1T_{1,l}=:\Lambda_1^{(l)}=\diag(i,-i,0).
\end{align*}
It fulfills the first order differential system
\begin{align*}
\widehat{W}^{(1,l)}_t+\left(\Lambda_0^{(l)}|\xi|^2+\Lambda_1^{(l)}|\xi|\right)\widehat{W}^{(1,l)}+A_0^{(1,l)}\widehat{W}^{(1,l)}=0,
\end{align*}
 where we employed
 \begin{align*}
 T_{1,l}^{-1}\Lambda_0^{(l)}T_{1,l}=T_{1,l}^{-1}T_{1,l}\Lambda_0^{(l)}=\Lambda_0^{(l)}
 \end{align*}
and
\begin{align*}
A_0^{(1,l)}:=T_{1,l}^{-1}A_0T_{1,l}=\frac{1}{2}\begin{pmatrix}
	0 & 0 & 1-i\\
	0 & 0 & 1+i\\
	-1-i & -1+i &0
\end{pmatrix}.
\end{align*}

Next, by taking the auxiliary matrix
\begin{align}\label{T2-l}
T_{2,l}:=I_3+N_{2,l}\ \ \mbox{with}\ \ N_{2,l}:=\frac{|\xi|^{-2}}{2}
\begin{pmatrix}
0 & 0 & 1-i\\
0 & 0 & 1+i\\
1+i & 1-i & 0
\end{pmatrix},
\end{align}
we can define the further unknown $\widehat{W}^{(2,l)}:=T_{2,l}^{-1}\widehat{W}^{(1,l)}$, which solves
\begin{align*}
\widehat{W}^{(2,l)}_t&+\left(\Lambda_0^{(l)}|\xi|^2+\Lambda_1^{(l)}|\xi|\right)\widehat{W}^{(2,l)}+T_{2,l}^{-1}\left(A_0^{(1,l)}-[N_{2,l},\Lambda_0^{(l)}]|\xi|^2\right)\widehat{W}^{(2,l)}\\
&-N_{2,l}T_{2,l}^{-1}\Lambda_1^{(l)}|\xi|\widehat{W}^{(2,l)}+T_{2,l}^{-1}\left(\Lambda_1^{(l)}|\xi|+A_0^{(1,l)}\right)N_{2,l}\widehat{W}^{(2,l)}=0.
\end{align*}
With the aid of
\begin{align*}
A_0^{(1,l)}-[N_{2,l},\Lambda_0^{(l)}]|\xi|^2=:\Lambda_2^{(l)}=0_{3\times 3},
\end{align*}
the coupled system is reduced to
\begin{align*}
\widehat{W}^{(2,l)}_t+\left(\Lambda_0^{(l)}|\xi|^2+\Lambda_1^{(l)}|\xi|\right)\widehat{W}^{(2,l)}+\left(A_0^{(2,l)}|\xi|^{-1}+R_{2,l}\right)\widehat{W}^{(2,l)}=0,
\end{align*}
where the coefficient matrix is
\begin{align*}
A_0^{(2,l)}|\xi|^{-1}:=[\Lambda_1^{(l)},N_{2,l}]|\xi|=\frac{1}{2}\begin{pmatrix}
	0 & 0 & 1+i\\
	0 & 0 & 1-i\\
	1-i & 1+i & 0
\end{pmatrix}|\xi|^{-1},
\end{align*}
and its remainder $R_{2,l}=O(|\xi|^{-2})$ is denoted by
\begin{align*}
R_{2,l}:=N_{2,l}^2T_{2,l}^{-1}\Lambda_1^{(l)}|\xi|-N_{2,l}T_{2,l}^{-1}\Lambda_1^{(l)}N_{2,l}|\xi|+T_{2,l}^{-1}A_0^{(1,l)}N_{2,l}.
\end{align*}

We next introduce
\begin{align}\label{T3-l}
T_{3,l}:=I_3+N_{3,l}\ \ \mbox{with}\ \ N_{3,l}:=\frac{|\xi|^{-3}}{2}\begin{pmatrix}
0 & 0 & 1+i\\
0 & 0 & 1-i\\
-1+i & -1-i & 0
\end{pmatrix}.
\end{align}
Thus, the new unknown $\widehat{W}^{(3,l)}:=T_{3,l}^{-1}\widehat{W}^{(2,l)}$ satisfies
\begin{align*}
\widehat{W}^{(3,l)}_t&+\left(\Lambda_0^{(l)}|\xi|^2+\Lambda_1^{(l)}|\xi|\right)\widehat{W}^{(3,l)}+T_{3,l}^{-1}R_{2,l}T_{3,l}\widehat{W}^{(3,l)}\\
&+\left(T_{3,l}^{-1}A_0^{(2,l)}N_{3,l}|\xi|^{-1}+T_{3,l}^{-1}\Lambda_1^{(l)}N_{3,l}|\xi|-N_{3,l}T_{3,l}^{-1}\Lambda_1^{(l)}|\xi|\right)\widehat{W}^{(3,l)}=0,
\end{align*}
where we used the fact that
\begin{align*}
A_0^{(2,l)}|\xi|^{-1}-[N_{3,l},\Lambda_0^{(l)}]|\xi|^2=:\Lambda_3^{(l)}|\xi|^{-1}=0_{3\times 3}.
\end{align*}
Furthermore, let us take the non-diagonal part with the size $|\xi|^{-2}$, that is
\begin{align*}
A_0^{(3,l)}|\xi|^{-2}:=A_0^{(1,l)}N_{2,l}+[\Lambda_1^{(l)},N_{3,l}]|\xi|=\frac{1}{2}
\begin{pmatrix}
1 & -i & -1+i\\
i & 1 & -1-i\\
1+i & 1-i & -2
\end{pmatrix}|\xi|^{-2},
\end{align*}
which yields the following first order  system:
\begin{align*}
\widehat{W}^{(3,l)}_t+\left(\Lambda_0^{(l)}|\xi|^2+\Lambda_1^{(l)}|\xi|\right)\widehat{W}^{(3,l)}+\left(A_0^{(3,l)}|\xi|^{-2}+R_{3,l}\right)\widehat{W}^{(3,l)}=0
\end{align*}
with the suitable remainder $R_{3,l}=O(|\xi|^{-3})$.

Eventually, benefiting from the non-zero diagonal part of $A_0^{(3,l)}$, we may take the auxiliary matrices
\begin{align}
T_{4,l}&:=I_3+N_{4,l}\ \ \mbox{with}\ \ N_{4,l}:=\frac{|\xi|^{-4}}{2}
\left(\begin{array}{cc;{2pt/2pt}c}
0 & 0 & -1+i\\
0 & 0 & -1-i\\ \hdashline[2pt/2pt]
-1-i & -1+i & 0
\end{array}\right),\label{T4-l}\\
T_{5,l}&:=I_3+N_{5,l}\ \ \mbox{with}\ \ N_{5,l}:=\frac{|\xi|^{-3}}{4}
\left(\begin{array}{cc;{2pt/2pt}c}
	0 & 1 & 0\\
	1 & 0 & 0\\ \hdashline[2pt/2pt]
	0 & 0 & 0
\end{array}\right).\label{T5-l}
\end{align}
Since the identity of the first two diagonal elements of $\Lambda_0^{(l)}$, we construct the matrix $N_{4,l}$ to treat four elements (locating in $1^{\mathrm{st}}$ and $2^{\mathrm{nd}}$ rows with $3^{\mathrm{rd}}$ column; $3^{\mathrm{rd}}$ row with $1^{\mathrm{st}}$ and $2^{\mathrm{nd}}$ columns) in $A_0^{(3,l)}$, precisely,
\begin{align*}
A_0^{(3,l)}|\xi|^{-2}-[N_{4,l},\Lambda_0^{(l)}]|\xi|^2=:A_0^{(4,l)}|\xi|^{-2}=\frac{1}{2} \left(\begin{array}{cc;{2pt/2pt}c}
1 & -i & 0\\
i & 1 & 0\\ \hdashline[2pt/2pt]
0 & 0  & -2
\end{array}\right)|\xi|^{-2}.
\end{align*}
The rest of two elements (locating in $1^{\mathrm{st}}$ row with $2^{\mathrm{nd}}$ column; $2^{\mathrm{nd}}$ row with $1^{\mathrm{st}}$ column) in $A_0^{(3,l)}$ or $A_0^{(4,l)}$ can be compensated by the other matrix $T_{5,l}$ in the following way:
\begin{align*}
A_0^{(4,l)}|\xi|^{-2}-[N_{5,l},\Lambda_1^{(l)}]|\xi|=:\Lambda_4^{(l)}|\xi|^{-2}=\frac{1}{2}\diag(1,1,-2)|\xi|^{-2}.
\end{align*}
In the above treatment, we applied the two steps diagonalization procedure (see, for example, \cite[Section 4]{Liu-Reissig=2014}). Consequently, we derive
\begin{align*}
\widehat{W}^{(4,l)}_t+\left(\Lambda_0^{(l)}|\xi|^2+\Lambda_1^{(l)}|\xi|+\Lambda_4^{(l)}|\xi|^{-2}\right)\widehat{W}^{(4,l)}+R_{4,l}\widehat{W}^{(4,l)}=0
\end{align*}
with the suitable remainder $R_{4,l}=O(|\xi|^{-3})$.

In short, by carrying out four steps of diagonalization procedure, we have derived the diagonal matrix in the last coupled system
\begin{align*}
	\Lambda_0^{(l)}|\xi|^2+\Lambda_1^{(l)}|\xi|+\Lambda_4^{(l)}|\xi|^{-2}=\diag\left(i|\xi|+\frac{|\xi|^{-2}}{2},-i|\xi|+\frac{|\xi|^{-2}}{2},|\xi|^2-|\xi|^{-2}\right).
\end{align*}
Thanks to the positive real parts as well as the pairwise distinct diagonal elements of the previous matrix and the higher order term $R_{4,l}$ to be the remainder, we in Proposition \ref{Prop-Diag-Large} successfully get an asymptotic representation of the energy term in the Fourier space for large frequencies from the general philosophy proposed in \cite[Chapter 2]{Jachmann=2008}.
\begin{prop}\label{Prop-Diag-Large}
	The characteristic roots $\lambda_{1,2,3}$ of the coefficient matrix $A_0+A_1|\xi|+A_2|\xi|^2$ in the Cauchy problem \eqref{Fourier-First-Order} behave for $\xi\in\ml{Z}_{\extt}(N_0)$ as
	\begin{align*}
		\lambda_{1,2}=\mp i|\xi|-\frac{|\xi|^{-2}}{2}+O(|\xi|^{-3}),\ \	\lambda_3=-|\xi|^2+|\xi|^{-2}+O(|\xi|^{-3}).
	\end{align*}
	Moreover, the solution for $\xi\in\ml{Z}_{\extt}(N_0)$  to the Cauchy problem \eqref{Fourier-First-Order} has the following representation: 
	\begin{align*}
		\chi_{\extt}(\xi)\widehat{W}=\chi_{\extt}(\xi)T_l\diag\left(\mathrm{e}^{\lambda_1 t},\mathrm{e}^{\lambda_2 t},\mathrm{e}^{\lambda_3 t}\right)T_l^{-1}\widehat{W}_0
	\end{align*}
	with the coefficient $T_l=O(1)$ such that
	\begin{align*}
		T_l:=T_{1,l}(I_3+N_{2,l})(I_3+N_{3,l})(I_3+N_{4,l})(I_3+N_{5,l}),
	\end{align*}
	where the matrices $T_{1,l}$ and $N_{k,l}$ with $k=2,\dots,5$ are defined in \eqref{T1-l}-\eqref{T5-l}, respectively.
\end{prop}

\subsection{Exponential decay stability for bounded frequencies}
\hspace{5mm}We investigate in this subsection  an exponential decay stability for the solution in the Fourier space for bounded frequencies $\xi\in\ml{Z}_{\bdd}(\varepsilon_0,N_0)$. From Propositions \ref{Prop-Diag-Small} and \ref{Prop-Diag-Large}, we already obtained 
\begin{align*}
\mathrm{Re}\,\lambda_{1,2,3}<0\ \ \mbox{for}\ \ \xi\in\ml{Z}_{\intt}(\varepsilon_0)\cup\ml{Z}_{\extt}(N_0).
\end{align*}
 Thus, if the characteristic roots are not pure imaginary, thanks to the continuity of these roots with respect to $|\xi|$, it is easy to claim that 
 \begin{align*}
 	\mathrm{Re}\,\lambda_{1,2,3}<0\ \ \mbox{for}\ \ \xi\in\ml{Z}_{\bdd}(\varepsilon_0,N_0)\ \ \Rightarrow \ \ \mbox{an exponential decay estimate}.
 \end{align*}

Let us prove the precondition in the above that there is no pure imaginary characteristic root of the matrix $A_0+A_1|\xi|+A_2|\xi|^2$ by a contradiction argument. That is to say, we assume that there is a pure imaginary root $ia$ with $a\in\mb{R}\backslash\{0\}$ for $\xi\in\ml{Z}_{\bdd}(\varepsilon_0,N_0)$. This non-trivial real number $a$ should satisfy the equality
\begin{align*}
0=\det\left(A_0+A_1|\xi|+A_2|\xi|^2-ia I_{3}\right)=ia\left(a^2-1-|\xi|^2\right)+|\xi|^2\left(|\xi|^2-a^2\right),
\end{align*}
in other words,
\begin{align*}
a^2=1+|\xi|^2\ \ \mbox{as well as}\ \ a^2=|\xi|^2.
\end{align*}
It cannot be achieved simultaneously, which gives a contradiction. Then, we have the next result.
\begin{prop}\label{Prop-Bdd}
The solution for $\xi\in\ml{Z}_{\bdd}(\varepsilon_0,N_0)$ to the Cauchy problem \eqref{Fourier-First-Order} satisfies an exponential decay estimate
\begin{align*}
\chi_{\bdd}(\xi)|\widehat{W}|\lesssim\chi_{\bdd}(\xi)\mathrm{e}^{-ct}|\widehat{W}_0|
\end{align*}
with a positive constant $c$.
\end{prop}
\subsection{Decay estimate and asymptotic profile: Proof of Theorem \ref{Thm-Energy-Term}}
\hspace{5mm}Summarizing Propositions \ref{Prop-Diag-Small}-\ref{Prop-Bdd}, particularly,
\begin{align*}
\chi_{\intt}(\xi)|\widehat{W}|\lesssim\chi_{\intt}(\xi)\mathrm{e}^{-c|\xi|^4t}|\widehat{W}_0|\ \ \mbox{as well as}\ \ \chi_{\extt}(\xi)|\widehat{W}|\lesssim\chi_{\extt}(\xi)\mathrm{e}^{-c\frac{1}{1+|\xi|^2}t}|\widehat{W}_0|,
\end{align*}
 we are able to derive the following pointwise estimate in the Fourier space:
\begin{align}\label{Pointwise-Fourier}
|\widehat{W}|\lesssim\mathrm{e}^{-c\frac{|\xi|^4}{(1+|\xi|^2)^3}t}|\widehat{W}_0|,
\end{align}
which exactly coincides with the one in \cite[Theorem 1.2]{Ueda-2018}. The estimate \eqref{Pointwise-Fourier} shows the strict dissipativity of \emph{Type (2,3)}, whose definition is referred to \cite[Defintion 1.0.2]{Mori-Thesis=2016}. This type of weak dissipative structure has been recently observed for several interesting systems, e.g. the Timoshenko system with the Fourier law \cite{Mori-Kawashima=2014} or the Cattaneo law \cite{Mori-Kawashima=2016} of heat conductions.
\begin{remark}\label{Remark-Point-wise-without-Ueda}
Although different energy terms are considered between ours and \cite{Ueda-2018}, the decay property of energy terms involving $u_x,u_t,v$ in general for the hyperbolic-parabolic coupled system \eqref{Linear-Hyper-Para-Syst} has been given already. This is the reason for the same decay property between \eqref{Pointwise-Fourier} and \cite[Estimate (1.5)]{Ueda-2018}. It is remarkable that the constraint condition (cf. \cite[Equality (1.4) and Theorem 1.2]{Ueda-2018}) is not any restriction because it works for $\nabla u_0$ in that paper.
\end{remark}

\subsubsection{Upper bound decay estimate of the regularity-loss type}
\hspace{5mm}We first derive the upper bound decay estimate for the energy term $W(t,\cdot)$ in the $\dot{H}^{\sigma}$ norm with $\sigma\geqslant0$. The Plancherel theorem and the Hausdorff-Young inequality yield
\begin{align}\label{Est-W}
\|W(t,\cdot)\|_{\dot{H}^{\sigma}}&\lesssim \left\|\chi_{\intt}(\xi)|\xi|^{\sigma}\mathrm{e}^{-c|\xi|^4t}\widehat{W}_0\right\|_{L^2}+\mathrm{e}^{-ct}\left\|\chi_{\bdd}(\xi)|\xi|^{\sigma}\widehat{W}_0\right\|_{L^2}+\left\|\chi_{\extt}(\xi)|\xi|^{\sigma}\mathrm{e}^{-c|\xi|^{-2}t}\widehat{W}_0\right\|_{L^2}\notag\\
&\lesssim\left\|\chi_{\intt}(\xi)|\xi|^{\sigma}\mathrm{e}^{-c|\xi|^4t}\right\|_{L^2}\|W_0\|_{L^1}+\mathrm{e}^{-ct}\|W_0\|_{L^2}+\left\|\chi_{\extt}(\xi)|\xi|^{-\ell}\mathrm{e}^{-c|\xi|^{-2}t}\right\|_{L^{\infty}}\|W_0\|_{H^{\sigma+\ell}}\notag\\
&\lesssim\left(\int_0^{\varepsilon_0}r^{2\sigma+n-1}\mathrm{e}^{-2cr^4t}\mathrm{d}r\right)^{1/2}\|W_0\|_{L^1}+\mathrm{e}^{-ct}\|W_0\|_{L^2}+\sup\limits_{r\geqslant N_0}\left(r^{-\ell}\mathrm{e}^{-cr^{-2}t}\right)\|W_0\|_{H^{\sigma+\ell}}\notag\\
&\lesssim (1+t)^{-\frac{n+2\sigma}{8}}\|W_0\|_{L^1}+(1+t)^{-\frac{\ell}{2}}\|W_0\|_{H^{\sigma+\ell}}
\end{align}
with $\sigma\geqslant0$ and $\ell\geqslant0$, where we used the fact $|\xi|^{\sigma+\ell}\approx\langle \xi\rangle^{\sigma+\ell}$ for $\xi\in\ml{Z}_{\extt}(N_0)$ and the estimates
\begin{align}
\int_0^{\varepsilon_0}r^{n+2\sigma-1}\mathrm{e}^{-2cr^4t}\mathrm{d}r&=t^{-\frac{2\sigma +n}{4}}\int_0^{\varepsilon_0}(r^4t)^{\frac{n+2\sigma-1}{4}}\mathrm{e}^{-2c(r^4t)}\mathrm{d}(r^4t)^{\frac{1}{4}}\lesssim 
\begin{cases}
1&\mbox{if}\ \ t\leqslant 1,\\
t^{-\frac{n+2\sigma}{4}}&\mbox{if}\ \ t>1,\\
\end{cases}\label{Technique-small}\\
\sup\limits_{r\geqslant N_0}\left(r^{-\ell}\mathrm{e}^{-cr^{-2}t}\right)&=t^{-\frac{\ell}{2}}\sup\limits_{r\geqslant N_0}\left((r^{-2}t)^{\frac{\ell}{2}}\mathrm{e}^{-c(r^{-2}t)}\right)\lesssim 
\begin{cases}
	1&\mbox{if}\ \ t\leqslant 1,\\
	t^{-\frac{\ell}{2}}&\mbox{if}\ \ t>1.
\end{cases}\label{Technique-large}
\end{align}
In other words, our desired upper bound estimate \eqref{Upper-Bound-Est} arises.

\subsubsection{Large time asymptotic profile}
\hspace{5mm}We next study a large time asymptotic profile of the energy term $W(t,\cdot)$ in the $\dot{H}^{\sigma}$ norm with some refined estimates (in the sense of improved decay rate and regularity). Let us introduce
\begin{align*}
\mu_1^{(s)}:=-|\xi|^4\ \ \mbox{and}\ \ \mu_{2,3}^{(s)}:=\mp i-\frac{1\pm i}{2}|\xi|^2+\frac{2\pm i}{4}|\xi|^4,
\end{align*}
moreover,
\begin{align*}
\mu_{1,2}^{(l)}:=\mp i|\xi|-\frac{|\xi|^{-2}}{2}\ \ \mbox{and}\ \ \mu_3^{(l)}:=-|\xi|^2+|\xi|^{-2},
\end{align*}
which are the corresponding principal parts of the characteristic roots $\lambda_{1,2,3}$ with suitable orders for $\xi\in\ml{Z}_{\intt}(\varepsilon_0)$ [resp. $\xi\in\ml{Z}_{\extt}(N_0)$] in Proposition \ref{Prop-Diag-Small} [resp. Proposition \ref{Prop-Diag-Large}]. 

For one thing, the solution to the reference system \eqref{Reference-System-I} in the Fourier space is uniquely given by
\begin{align}\label{Vs}
\widehat{V}^{(s)}=\diag\left(\mathrm{e}^{\mu_1^{(s)}t},\mathrm{e}^{\mu_2^{(s)}t},\mathrm{e}^{\mu_3^{(s)}t}\right)\widehat{Q}_s^{-1}\widehat{W}_0.
\end{align}
Notice that $T_s=\widehat{Q}_s(I_3+N_{5,s})$. We may rewrite the solution in Proposition \ref{Prop-Diag-Small} by
\begin{align*}
\chi_{\intt}(\xi)\widehat{W}&=\chi_{\intt}(\xi)\widehat{Q}_s(I_3+N_{5,s})\diag\left(\mathrm{e}^{\lambda_1t},\mathrm{e}^{\lambda_2t},\mathrm{e}^{\lambda_3t}\right)(I_3+N_{5,s})^{-1}\widehat{Q}_s^{-1}\widehat{W}_0\\
&=\chi_{\intt}(\xi)\widehat{Q}_s\diag\left(\mathrm{e}^{\lambda_1t},\mathrm{e}^{\lambda_2t},\mathrm{e}^{\lambda_3t}\right)\widehat{Q}_s^{-1}\widehat{W}_0+\chi_{\intt}(\xi)\widehat{Q}_sN_{5,s}\diag\left(\mathrm{e}^{\lambda_1t},\mathrm{e}^{\lambda_2t},\mathrm{e}^{\lambda_3t}\right)T_s^{-1}\widehat{W}_0\\
&\quad-\chi_{\intt}(\xi)\widehat{Q}_s\diag\left(\mathrm{e}^{\lambda_1t},\mathrm{e}^{\lambda_2t},\mathrm{e}^{\lambda_3t}\right)N_{5,s}T_s^{-1}\widehat{W}_0,
\end{align*}
where we used the next identity again:
\begin{align*}
(I_3+N_{5,s})^{-1}=I_3-N_{5,s}(I_3+N_{5,s})^{-1}.
\end{align*}
Let us recall $\widehat{Q}_s=O(1)$, $T_s=O(1)$, $N_{5,s}=O(|\xi|^4)$ and $\lambda_k-\mu_k^{(s)}=O(|\xi|^5)$ with $k=1,2,3$ for $\xi\in\ml{Z}_{\intt}(\varepsilon_0)$. A direct subtraction shows
\begin{align}\label{Error-EST-small}
\chi_{\intt}(\xi)\left|\widehat{W}-\widehat{Q}_s\widehat{V}^{(s)}\right|&\lesssim\chi_{\intt}(\xi)\max\limits_{k=1,2,3}\left|\mathrm{e}^{\lambda_kt}-\mathrm{e}^{\mu_k^{(s)}t}\right||\widehat{W}_0|+\chi_{\intt}(\xi)|\xi|^4\max\limits_{k=1,2,3}\mathrm{e}^{\lambda_kt}|\widehat{W}_0|\notag\\
&\lesssim\chi_{\intt}(\xi)\max\limits_{k=1,2,3}\left|(\lambda_k-\mu_k^{(s)})t\mathrm{e}^{\mu_k^{(s)}t}\int_0^1\mathrm{e}^{(\lambda_k-\mu_k^{(s)})t\eta}\mathrm{d}\eta\right|+\chi_{\intt}(\xi)|\xi|^4\max\limits_{k=1,2,3}\mathrm{e}^{\lambda_kt}|\widehat{W}_0|\notag\\
&\lesssim\chi_{\intt}(\xi)\left(|\xi|^5t+|\xi|^4\right)\mathrm{e}^{-c|\xi|^4t}|\widehat{W}_0|\notag\\
&\lesssim\chi_{\intt}(\xi)|\xi|\mathrm{e}^{-c|\xi|^4t}|\widehat{W}_0|.
\end{align}
It leads to 
\begin{align}\label{Error-small}
\left\|\chi_{\intt}(D)\left(W(t,\cdot)-Q_sV^{(s)}(t,\cdot)\right)\right\|_{\dot{H}^{\sigma}}&\lesssim\left\|\chi_{\intt}(\xi)|\xi|^{\sigma+1}\mathrm{e}^{-c|\xi|^4t}\right\|_{L^2}\|W_0\|_{L^1}\notag\\
&\lesssim (1+t)^{-\frac{n+2\sigma}{8}-\frac{1}{4}}\|W_0\|_{L^1},
\end{align}
where the estimate \eqref{Technique-small} with $\sigma\to \sigma+1$ was applied.

For another, the solution to the reference system \eqref{Reference-System-II} in the Fourier space is represented by
\begin{align*}
\widehat{V}^{(l)}=\diag\left(\mathrm{e}^{\mu_1^{(l)}t},\mathrm{e}^{\mu_2^{(l)}t},\mathrm{e}^{\mu_3^{(l)}t}\right)\widehat{Q}_l^{-1}\widehat{W}_0.
\end{align*}
Similarly to the case for small frequencies, we are able to deduce
\begin{align*}
\chi_{\extt}(\xi)\left|\widehat{W}-\widehat{Q}_l\widehat{V}^{(l)}\right|&\lesssim\chi_{\extt}(\xi)\left(|\xi|^{-3}t+|\xi|^{-3}\right)\mathrm{e}^{-c|\xi|^{-2}t}|\widehat{W}_0|\\
&\lesssim\chi_{\extt}(\xi)|\xi|^{-1}\mathrm{e}^{-c|\xi|^{-2}t}|\widehat{W}_0|,
\end{align*}
where we mainly used $N_{5,l}=O(|\xi|^{-3})$ and $\lambda_k-\mu_k^{(l)}=O(|\xi|^{-3})$ with $k=1,2,3$ for $\xi\in\ml{Z}_{\extt}(N_0)$. That is to say, 
\begin{align}\label{Error-large}
\left\|\chi_{\extt}(D)\left(W(t,\cdot)-Q_lV^{(l)}(t,\cdot)\right)\right\|_{\dot{H}^{\sigma}}&\lesssim\left\|\chi_{\extt}(\xi)|\xi|^{-\ell-\frac{1}{2}}\mathrm{e}^{-c|\xi|^{-2}t}\right\|_{L^{\infty}}\|W_0\|_{H^{\sigma+\ell-\frac{1}{2}}}\notag\\
&\lesssim (1+t)^{-\frac{\ell}{2}-\frac{1}{4}}\|W_0\|_{H^{\sigma+\ell-\frac{1}{2}}},
\end{align}
where the estimate \eqref{Technique-large} with $\ell\to\ell+\frac{1}{2}$ was applied.

Summarizing the derived estimates \eqref{Error-small} and \eqref{Error-large} and the exponential decay estimate for bounded frequencies, one can prove our desired error estimate \eqref{Error-ALL} immediately. 

\subsubsection{Optimal decay estimate for large time}
\hspace{5mm}Let us additionally consider in this part the suitable higher regularity $\ell>\frac{n}{4}+\frac{\sigma}{2}$ with $n\geqslant 1$ and $\sigma\geqslant 0$. From the derived estimate \eqref{Upper-Bound-Est}, we may get
\begin{align}\label{Upper-02}
\|W(t,\cdot)\|_{\dot{H}^{\sigma}}\lesssim t^{-\frac{n+2\sigma}{8}}\|W_0\|_{H^{\sigma+\ell}\cap L^1}
\end{align}
for large time $t\gg1$. We next demonstrate the optimality of the last decay rate $t^{-\frac{n+2\sigma}{8}}$, that is the sharp lower bound estimate with the $L^1$ data. 

Recalling the representation \eqref{Vs} and the differential operator $Q_s=Q_s(|D|)$, we are able to rewrite the profile by
\begin{align*}
\chi_{\intt}(D)|D|^{\sigma}Q_sV^{(s)}(t,x)=:\ml{\mb{L}}_{0,\sigma}(t,|D|)W_0(x):=\ml{L}_{0,\sigma}(t,x)\ast_{(x)} W_0(x).
\end{align*}
According to the mean value theorem, there exists $\theta_0\in(0,1)$ such that
\begin{align*}
|\ml{L}_{0,\sigma}(t,x-y)-\ml{L}_{0,\sigma}(t,x)|\lesssim|y|\,|\nabla\ml{L}_{0,\sigma}(t,x-\theta_0y)|.
\end{align*}
We now separate the domain of integral into $\{|y|\leqslant t^{\frac{1}{8}}\}$ and $\{|y|\geqslant t^{\frac{1}{8}}\}$, precisely, 
\begin{align*}
&\|\ml{\mb{L}}_{0,\sigma}(t,|D|)W_0(\cdot)-\ml{L}_{0,\sigma}(t,\cdot)P_{W_0}\|_{L^2}\\
&=\left\|\ml{L}_{0,\sigma}(t,\cdot)\ast_{(x)}W_0(\cdot)-\ml{L}_{0,\sigma}(t,\cdot)P_{W_0}\right\|_{L^2}\\
&=\left\|\int_{\mb{R}^n}\ml{L}_{0,\sigma}(t,\cdot-y)W_0(y)\mathrm{d}y-\ml{L}_{0,\sigma}(t,\cdot)\int_{\mb{R}^n}W_0(y)\mathrm{d}y\right\|_{L^2}\\
&\leqslant\left\|\int_{|y|\leqslant t^{\frac{1}{8}}}|\ml{L}_{0,\sigma}(t,\cdot-y)-\ml{L}_{0,\sigma}(t,\cdot)|\,|W_0(y)|\mathrm{d}y\right\|_{L^2}+\left\|\int_{|y|\geqslant t^{\frac{1}{8}}}\big(|\ml{L}_{0,\sigma}(t,\cdot-y)|+|\ml{L}_{0,\sigma}(t,\cdot)|\big)|W_0(y)|\mathrm{d}y \right\|_{L^2}.
\end{align*}
As a consequence,
\begin{align}
&\|\ml{\mb{L}}_{0,\sigma}(t,|D|)W_0(\cdot)-\ml{L}_{0,\sigma}(t,\cdot)P_{W_0}\|_{L^2}\notag\\
&\lesssim \left\|\int_{|y|\leqslant t^{\frac{1}{8}}}|y|\,|\nabla\ml{L}_{0,\sigma}(t,\cdot-\theta_0y)|\,|W_0(y)|\mathrm{d}y\right\|_{L^2}+\|\ml{L}_{0,\sigma}(t,\cdot)\|_{L^2}\int_{|y|\geqslant t^{\frac{1}{8}}}|W_0(y)|\mathrm{d}y\notag\\
&\lesssim t^{\frac{1}{8}}\|\ml{L}_{0,\sigma+1}(t,\cdot)\|_{L^2}\|W_0\|_{L^1}+\|\ml{L}_{0,\sigma}(t,\cdot)\|_{L^2}\int_{|y|\geqslant t^{\frac{1}{8}}}|W_0(y)|\mathrm{d}y\notag\\
&\lesssim t^{\frac{1}{8}-\frac{n+2(\sigma+1)}{8}}\|W_0\|_{L^1}+o(t^{-\frac{n+2\sigma}{8}})\notag\\
&=o(t^{-\frac{n+2\sigma}{8}})\label{Est-Error-01}
\end{align}
for large time $t\gg1$, where we used
\begin{align*}
\|\ml{L}_{0,\sigma}(t,\cdot)\|_{L^2}^2\simeq \left\|\chi_{\intt}(\xi)|\xi|^{\sigma}\mathrm{e}^{-c|\xi|^4t}\right\|_{L^2}^2\simeq t^{-\frac{n+2\sigma}{4}}
\end{align*}
and the $L^1$ integrability of $W_0$ leading to
\begin{align*}
\int_{|y|\geqslant t^{\frac{1}{8}}}|W_0(y)|\mathrm{d}y=o(1)\ \ \mbox{when}\ \ t\gg1.
\end{align*}
Moreover, a combination of \eqref{Est-W} and \eqref{Error-EST-small} shows
\begin{align}
\||D|^{\sigma}W(t,\cdot)-\ml{\mb{L}}_{0,\sigma}(t,|D|)W_0(\cdot)\|_{L^2}&\lesssim t^{-\frac{n+2\sigma}{8}}\left(t^{-\frac{1}{4}}\|W_0\|_{L^1}+t^{\frac{n+2\sigma}{8}-\frac{\ell}{2}}\|W_0\|_{H^{\sigma+\ell}}\right)\notag\\
&=o(t^{-\frac{n+2\sigma}{8}})\label{Est-Error-02}
\end{align}
for large time $t\gg1$, thanks to $\frac{n+2\sigma}{8}<\frac{\ell}{2}$ in our consideration. Applying the triangle inequality with \eqref{Est-Error-01} and \eqref{Est-Error-02}, we claim the following large time behavior:
\begin{align*}
\||D|^{\sigma}W(t,\cdot)-\ml{L}_{0,\sigma}(t,\cdot)P_{W_0}\|_{L^2}=o(t^{-\frac{n+2\sigma}{8}}).
\end{align*}
Finally, we employ the triangle inequality again to see
\begin{align}\label{Est-07}
\|W(t,\cdot)\|_{\dot{H}^{\sigma}}&\geqslant \|\ml{L}_{0,\sigma}(t,\cdot)\|_{L^2}|P_{W_0}|-\||D|^{\sigma}W(t,\cdot)-\ml{L}_{0,\sigma}(t,\cdot)P_{W_0}\|_{L^2}\notag\\
&\gtrsim t^{-\frac{n+2\sigma}{8}}|P_{W_0}|-o(t^{-\frac{n+2\sigma}{8}})\notag\\
&\gtrsim t^{-\frac{n+2\sigma}{8}}|P_{W_0}|
\end{align}
for large time $t\gg1$. In the view of the upper bound estimate \eqref{Upper-02}, it completes our proof of sharp lower bound estimate.

\section{Large time asymptotic behavior for the solution}\setcounter{equation}{0}\label{Section-Solution}
\hspace{5mm}Let us act the heat operator $\partial_t-\Delta$ on the equation \eqref{Linear-Hyper-Para-Syst}$_1$, and apply the equation \eqref{Linear-Hyper-Para-Syst}$_2$ to eliminate the unknown $v$. Then, we arrive at
\begin{align}\label{New-Model}
(\partial_t-\Delta)(u_{tt}-\Delta u)+u_t=0.
\end{align}
To the knowledge of the authors, the third order (in time) PDE \eqref{New-Model} does not been studied. This model can be understood by the heat operator $\partial_t-\Delta$ acting on the classical wave operator $\partial_t^2-\Delta$ with an additional friction term $u_t$. However, different from the classical damped wave equation, this term $u_t$ here posses the weak dissipative structure in the regularity-loss type only. It seems interesting to understand some influence of this weak dissipation on qualitative properties of the solution $u$. In other words, the unknown $u$ solves the  Cauchy problem for the following third order (in time) evolution equation:
\begin{align}\label{Third-Order-PDE}
\begin{cases}
u_{ttt}-\Delta u_{tt}-\Delta u_t+u_t+\Delta^2u=0,&x\in\mb{R}^n,\ t>0,\\
u(0,x)=u_0(x),\ u_t(0,x)=u_1(x),\ u_{tt}(0,x)=u_2(x),&x\in\mb{R}^n,
\end{cases}
\end{align}
where the third Cauchy data is given by \eqref{Linear-Hyper-Para-Syst}$_1$ with $t=0$ such that
\begin{align*}
u_2(x):=\Delta u_0(x)-v_0(x).
\end{align*}
Our aim in this section is to investigate the higher order PDE \eqref{Third-Order-PDE} by using the WKB analysis and the Fourier analysis.
\subsection{Representation of the solution in the Fourier space}
\hspace{5mm}A direct application of the partial Fourier transform with respect to the spatial variable, i.e. $\widehat{u}=\widehat{u}(t,\xi)$ such that $\widehat{u}=\ml{F}_{x\to\xi}(u)$, for the Cauchy problem \eqref{Third-Order-PDE} implies
\begin{align*}
	\begin{cases}
		\widehat{u}_{ttt}+|\xi|^2\widehat{u}_{tt}+(1+|\xi|^2)\widehat{u}_t+|\xi|^4\widehat{u}=0,&\xi\in\mb{R}^n,\ t>0,\\
		\widehat{u}(0,\xi)=\widehat{u}_0(\xi),\ \widehat{u}_t(0,\xi)=\widehat{u}_1(\xi),\ \widehat{u}_{tt}(0,\xi)=\widehat{u}_2(\xi),&\xi\in\mb{R}^n,
	\end{cases}
\end{align*}
whose characteristic equation with its roots $\lambda_{1,2,3}$ is addressed by
\begin{align}\label{Cubic-Eq}
\lambda^3+|\xi|^2\lambda^2+(1+|\xi|^2)\lambda+|\xi|^4=0.
\end{align}
It seems a challenging work to analyze the precise roots of this $|\xi|$-dependent cubic equation by Cardano's formula. Thus, we may use some asymptotic expansions for $\xi\in\ml{Z}_{\intt}(\varepsilon_0)\cup\ml{Z}_{\extt}(N_0)$ and a stable analysis for $\xi\in\ml{Z}_{\bdd}(\varepsilon_0,N_0)$ later.

We now apply the Taylor expansion  when $\xi\in\ml{Z}_{\intt}(\varepsilon_0)$ for the characteristic roots, that is
\begin{align*}
\lambda_k=\sum\limits_{j=0}^{+\infty}C_{k,j}^{(s)}|\xi|^j\ \ \mbox{with}\ \ C_{k,j}^{(s)}\in\mb{C}\ \ \mbox{for}\ \ k=1,2,3,
\end{align*}
to the cubic equation \eqref{Cubic-Eq} with analytic coefficients with respect to $|\xi|$. By straightforward but tedious computations, we have the asymptotic expansions
\begin{align*}
\lambda_1=-|\xi|^4+O(|\xi|^6),\ \ \lambda_{2,3}=\mp i-\frac{1\pm i}{2}|\xi|^2+O(|\xi|^4),
\end{align*}
which exactly coincide with those in Proposition \ref{Prop-Diag-Small}. Moreover, the discriminant of the cubic \eqref{Cubic-Eq} is $\triangle_{\mathrm{Dis}}=-4+O(|\xi|^2)<0$ when $\xi\in\ml{Z}_{\intt}(\varepsilon_0)$. It means that $\lambda_{2,3}$ are two complex conjugate roots, which can be rewritten by $\lambda_{2,3}=\lambda_{\mathrm{R}}\pm i\lambda_{\mathrm{I}}$ with
\begin{align*}
\lambda_{\mathrm{R}}=-\frac{|\xi|^2}{2}+O(|\xi|^4)\ \ \mbox{and}\ \ \lambda_{\mathrm{I}}=-1-\frac{|\xi|^2}{2}+O(|\xi|^4).
\end{align*}
Recalling the general solution's formula for third order (in time) PDEs with a pair of conjugate characteristic roots deduced in \cite[Equation (2.3)]{Chen-Takeda=2023}, the solution for small frequencies is explicit represented by
\begin{align*}
	\chi_{\intt}(\xi)\widehat{u}
	&=\chi_{\intt}(\xi)\left[\frac{(|\xi|^2-\lambda_{\mathrm{I}}^2-\lambda_{\mathrm{R}}^2)\widehat{u}_0+2\lambda_{\mathrm{R}}\widehat{u}_1+\widehat{v}_0}{2\lambda_{\mathrm{R}}\lambda_1-\lambda_{\mathrm{I}}^2-\lambda_{\mathrm{R}}^2-\lambda_1^2}\mathrm{e}^{\lambda_1t}\right.\\
	&\qquad\qquad\quad+\frac{(2\lambda_{\mathrm{R}}\lambda_1-\lambda_1^2-|\xi|^2)\widehat{u}_0-2\lambda_{\mathrm{R}}\widehat{u}_1-\widehat{v}_0}{2\lambda_{\mathrm{R}}\lambda_1-\lambda_{\mathrm{I}}^2-\lambda_{\mathrm{R}}^2-\lambda_1^2}\cos(\lambda_{\mathrm{I}}t)\mathrm{e}^{\lambda_{\mathrm{R}}t}\\
	&\qquad\qquad\quad+\frac{[\lambda_1(\lambda_{\mathrm{R}}\lambda_1+\lambda_{\mathrm{I}}^2-\lambda_{\mathrm{R}}^2)+|\xi|^2(\lambda_{\mathrm{R}}-\lambda_1)]\widehat{u}_0}{\lambda_{\mathrm{I}}(2\lambda_{\mathrm{R}}\lambda_1-\lambda_{\mathrm{I}}^2-\lambda_{\mathrm{R}}^2-\lambda_1^2)}\sin(\lambda_{\mathrm{I}}t)\mathrm{e}^{\lambda_{\mathrm{R}}t}\\
	&\qquad\qquad\quad\left.+\frac{(\lambda_{\mathrm{R}}^2-\lambda_{\mathrm{I}}^2-\lambda_1^2)\widehat{u}_1+(\lambda_{\mathrm{R}}-\lambda_1)\widehat{v}_0}{\lambda_{\mathrm{I}}(2\lambda_{\mathrm{R}}\lambda_1-\lambda_{\mathrm{I}}^2-\lambda_{\mathrm{R}}^2-\lambda_1^2)}\sin(\lambda_{\mathrm{I}}t)\mathrm{e}^{\lambda_{\mathrm{R}}t}\right],
\end{align*}
where we used the relation $\widehat{u}_2=-|\xi|^2\widehat{u}_0-\widehat{v}_0$.

Concerning $\xi\in\ml{Z}_{\extt}(N_0)$, we may employ the Laurent expansion, because of the smallness of $|\xi|^{-1}$, for the characteristic roots, namely,
\begin{align*}
\lambda_k=\sum\limits_{j=-2}^{+\infty} C_{k,j}^{(l)}|\xi|^{-j}\ \ \mbox{with}\ \ C_{k,j}^{(l)}\in\mb{C}\ \ \mbox{for}\ \ k=1,2,3,
\end{align*}
to the cubic equation \eqref{Cubic-Eq}. Similarly to the case for small frequencies, we immediately derive the asymptotic expansions
\begin{align*}
\lambda_{1,2}=\mp i|\xi|-\frac{|\xi|^{-2}}{2}+O(|\xi|^{-3}),\ \ \lambda_3=-|\xi|^2+O(|\xi|^{-2}),
\end{align*}
which exactly coincide with those in Proposition \ref{Prop-Diag-Large}. In this situation, the discriminant of the cubic \eqref{Cubic-Eq} is $\triangle_{\mathrm{Dis}}=-4|\xi|^{10}+O(|\xi|^8)<0$ when $\xi\in\ml{Z}_{\extt}(N_0)$. That is to say, $\lambda_{1,2}$ are two complex conjugate roots, which can be rewritten by $\lambda_{1,2}=\lambda_{\mathrm{R}}\pm i\lambda_{\mathrm{I}}$ with
\begin{align*}
\lambda_{\mathrm{R}}=-\frac{|\xi|^{-2}}{2}+O(|\xi|^{-3})\ \ \mbox{and}\ \ \lambda_{\mathrm{I}}=-|\xi|+O(|\xi|^{-3}).
\end{align*}
Therefore, the solution is expressed by
\begin{align*}
	\chi_{\extt}(\xi)\widehat{u}
	&=\chi_{\extt}(\xi)\left[\frac{(|\xi|^2-\lambda_{\mathrm{I}}^2-\lambda_{\mathrm{R}}^2)\widehat{u}_0+2\lambda_{\mathrm{R}}\widehat{u}_1+\widehat{v}_0}{2\lambda_{\mathrm{R}}\lambda_3-\lambda_{\mathrm{I}}^2-\lambda_{\mathrm{R}}^2-\lambda_3^2}\mathrm{e}^{\lambda_3t}\right.\\
	&\qquad\qquad\quad+\frac{(2\lambda_{\mathrm{R}}\lambda_3-\lambda_3^2-|\xi|^2)\widehat{u}_0-2\lambda_{\mathrm{R}}\widehat{u}_1-\widehat{v}_0}{2\lambda_{\mathrm{R}}\lambda_3-\lambda_{\mathrm{I}}^2-\lambda_{\mathrm{R}}^2-\lambda_3^2}\cos(\lambda_{\mathrm{I}}t)\mathrm{e}^{\lambda_{\mathrm{R}}t}\\
	&\qquad\qquad\quad+\frac{[\lambda_3(\lambda_{\mathrm{R}}\lambda_3+\lambda_{\mathrm{I}}^2-\lambda_{\mathrm{R}}^2)+|\xi|^2(\lambda_{\mathrm{R}}-\lambda_3)]\widehat{u}_0}{\lambda_{\mathrm{I}}(2\lambda_{\mathrm{R}}\lambda_3-\lambda_{\mathrm{I}}^2-\lambda_{\mathrm{R}}^2-\lambda_3^2)}\sin(\lambda_{\mathrm{I}}t)\mathrm{e}^{\lambda_{\mathrm{R}}t}\\
	&\qquad\qquad\quad\left.+\frac{(\lambda_{\mathrm{R}}^2-\lambda_{\mathrm{I}}^2-\lambda_3^2)\widehat{u}_1+(\lambda_{\mathrm{R}}-\lambda_3)\widehat{v}_0}{\lambda_{\mathrm{I}}(2\lambda_{\mathrm{R}}\lambda_3-\lambda_{\mathrm{I}}^2-\lambda_{\mathrm{R}}^2-\lambda_3^2)}\sin(\lambda_{\mathrm{I}}t)\mathrm{e}^{\lambda_{\mathrm{R}}t}\right].
\end{align*}

Lastly for bounded frequencies $\xi\in\ml{Z}_{\bdd}(\varepsilon_0,N_0)$, the exponentially decay estimate derived in Proposition \ref{Prop-Bdd} yields
\begin{align}\label{Est-u-BDD}
\chi_{\bdd}(\xi)|\widehat{u}|\lesssim\chi_{\bdd}(\xi)\mathrm{e}^{-ct}(|\widehat{u}_0|+|\widehat{u}_1|+|\widehat{v}_0|),
\end{align}
because of $\varepsilon_0\leqslant|\xi|\leqslant N_0$.

\subsection{Refined estimates in the Fourier space}
\hspace{5mm}According to the asymptotic expansions for small frequencies, we notice that the leading term of $\chi_{\intt}(\xi)\widehat{u}$ is given by
\begin{align*}
\chi_{\intt}(\xi)\widehat{J}_s:=\chi_{\intt}(\xi)\frac{-\lambda_{\mathrm{I}}^2\widehat{u}_0+\widehat{v}_0}{2\lambda_{\mathrm{R}}\lambda_1-\lambda_{\mathrm{I}}^2-\lambda_{\mathrm{R}}^2-\lambda_1^2}\mathrm{e}^{\lambda_1t},
\end{align*}
due to the error estimate
\begin{align}\label{Error-01}
\chi_{\intt}(\xi)|\widehat{u}-\widehat{J}_s|&\lesssim\chi_{\intt}(\xi)\left(|\xi|^2\mathrm{e}^{-c|\xi|^4t}(|\widehat{u}_0|+|\widehat{u}_1|)+\mathrm{e}^{-c|\xi|^2t}(|\xi|^2|\widehat{u}_0|+|\widehat{u}_1|+|\widehat{v}_0|)\right)\notag\\
&\lesssim\chi_{\intt}(\xi)\left(|\xi|^2\mathrm{e}^{-c|\xi|^4t}+\mathrm{e}^{-c|\xi|^2t}\right)(|\widehat{u}_0|+|\widehat{u}_1|+|\widehat{v}_0|),
\end{align}
and the sharp estimate
\begin{align*}
\chi_{\intt}(\xi)|\widehat{J}_s|\lesssim\chi_{\intt}(\xi)\mathrm{e}^{-c|\xi|^4t}(|\widehat{u}_0|+|\widehat{v}_0|).
\end{align*}
It implies from the triangle inequality
\begin{align}\label{Est-u-INT}
\chi_{\intt}(\xi)|\widehat{u}|\leqslant\chi_{\intt}(\xi)\left(|\widehat{u}-\widehat{J}_s|+|\widehat{J}_s|\right)\lesssim \chi_{\intt}(\xi)\mathrm{e}^{-c|\xi|^4t}(|\widehat{u}_0|+|\widehat{u}_1|+|\widehat{v}_0|).
\end{align}
Motivated by \eqref{Error-EST-small} for treating with the difference between two exponentially decay functions, we may find the following approximation:
\begin{align}
&\chi_{\intt}(\xi)\left|\widehat{J}_s-\mathrm{e}^{-|\xi|^4t}(\widehat{u}_0-\widehat{v}_0)\right|\notag\\
&\leqslant\chi_{\intt}(\xi)\left|\frac{-\lambda_{\mathrm{I}}^2\widehat{u}_0+\widehat{v}_0}{2\lambda_{\mathrm{R}}\lambda_1-\lambda_{\mathrm{I}}^2-\lambda_{\mathrm{R}}^2-\lambda_1^2}-(\widehat{u}_0-\widehat{v}_0) \right|\mathrm{e}^{\lambda_1t}+\chi_{\intt}(\xi)\left|\mathrm{e}^{\lambda_1t}-\mathrm{e}^{-|\xi|^4t}\right||\widehat{u}_0-\widehat{v}_0|\notag\\
&\lesssim\chi_{\intt}(\xi)\left(|\xi|^2\mathrm{e}^{-c|\xi|^4t}(|\xi|^2|\widehat{u}_0|+|\widehat{v}_0|)+|\xi|^6t\mathrm{e}^{-c|\xi|^4t}(|\widehat{u}_0|+|\widehat{v}_0|)\right)\notag\\
&\lesssim\chi_{\intt}(\xi)|\xi|^2\mathrm{e}^{-c|\xi|^4t}(|\widehat{u}_0|+|\widehat{v}_0|),\label{Error-02}
\end{align}
where we considered the refined relations $1-\lambda_{\mathrm{I}}^2=O(|\xi|^2)$ and $\lambda_1+|\xi|^4=O(|\xi|^6)$ for $\xi\in\ml{Z}_{\intt}(\varepsilon_0)$. Hence, the combination of \eqref{Error-01} and \eqref{Error-02} yields
\begin{align*}
\chi_{\intt}(\xi)\left|\widehat{u}-\mathrm{e}^{-|\xi|^4t}(\widehat{u}_0-\widehat{v}_0)\right|\lesssim\chi_{\intt}(\xi)\left(|\xi|^2\mathrm{e}^{-c|\xi|^4t}+\mathrm{e}^{-c|\xi|^2t}\right)(|\widehat{u}_0|+|\widehat{u}_1|+|\widehat{v}_0|).
\end{align*}

Let us turn to the case for large frequencies. The leading term of $\chi_{\extt}(\xi)\widehat{u}$ is denoted by
\begin{align*}
\chi_{\extt}(\xi)\widehat{J}_l:=\chi_{\extt}(\xi)\frac{-\lambda_3^2\widehat{u}_0}{2\lambda_{\mathrm{R}}\lambda_3-\lambda_{\mathrm{I}}^2-\lambda_{\mathrm{R}}^2-\lambda_3^2}\cos(\lambda_{\mathrm{I}}t)\mathrm{e}^{\lambda_{\mathrm{R}}t},
\end{align*}
due to the error estimate
\begin{align}\label{Error-03}
\chi_{\extt}(\xi)|\widehat{u}-\widehat{J}_l|&\lesssim\chi_{\extt}(\xi)\left(|\xi|^{-4}\mathrm{e}^{-c|\xi|^2t}(|\xi|^{-2}|\widehat{u}_0|+|\xi|^{-2}|\widehat{u}_1|+|\widehat{v}_0|)\right.\notag\\
&\qquad\qquad\quad\left.+|\xi|^{-1}\mathrm{e}^{-c|\xi|^{-2}t}(|\xi|^{-2}|\widehat{u}_0|+|\widehat{u}_1|+|\xi|^{-2}|\widehat{v}_0|)\right)\notag\\
&\lesssim\chi_{\extt}(\xi)|\xi|^{-1}\mathrm{e}^{-c|\xi|^{-2}t}(|\xi|^{-2}|\widehat{u}_0|+|\widehat{u}_1|+|\xi|^{-2}|\widehat{v}_0|),
\end{align}
where we used the crucial cancellation $|\xi|^2-\lambda_{\mathrm{I}}^2=O(|\xi|^{-2})$ twice for $\xi\in\ml{Z}_{\extt}(N_0)$  to get the coefficient of $\widehat{u}_0$, and the sharp estimate
\begin{align*}
\chi_{\extt}(\xi)|\widehat{J}_l|\lesssim\chi_{\extt}(\xi)\mathrm{e}^{-c|\xi|^{-2}t}|\widehat{u}_0|.
\end{align*}
It implies from the triangle inequality
\begin{align}\label{Est-u-EXT}
	\chi_{\extt}(\xi)|\widehat{u}|\leqslant\chi_{\extt}(\xi)\left(|\widehat{u}-\widehat{J}_l|+|\widehat{J}_l|\right)\lesssim \chi_{\extt}(\xi)\mathrm{e}^{-c|\xi|^{-2}t}(|\widehat{u}_0|+|\xi|^{-1}|\widehat{u}_1|+|\xi|^{-3}|\widehat{v}_0|).
\end{align}
Then, one may derive
\begin{align}\label{Error-04}
&\chi_{\extt}(\xi)\left|\widehat{J}_l-\cos(|\xi|t)\mathrm{e}^{-\frac{1}{2}|\xi|^{-2}t}\widehat{u}_0\right|\notag\\
&\leqslant\chi_{\extt}(\xi)\left|\frac{-\lambda_3^2\widehat{u}_0}{2\lambda_{\mathrm{R}}\lambda_3-\lambda_{\mathrm{I}}^2-\lambda_{\mathrm{R}}^2-\lambda_3^2}-\widehat{u}_0\right||\cos(\lambda_{\mathrm{I}}t)|\mathrm{e}^{\lambda_{\mathrm{R}}t}+\chi_{\extt}(\xi)\big|\cos(\lambda_{\mathrm{I}}t)-\cos(-|\xi|t)\big|\mathrm{e}^{\lambda_{\mathrm{R}}t}|\widehat{u}_0|\notag\\
&\quad+\chi_{\extt}(\xi)|\cos(|\xi|t)|\left|\mathrm{e}^{\lambda_{\mathrm{R}}t}-\mathrm{e}^{-\frac{1}{2}|\xi|^{-2}t}\right||\widehat{u}_0|\notag\\
&\lesssim\chi_{\extt}(\xi)\left(|\xi|^{-2}+|\xi|^{-3}t\right)\mathrm{e}^{-c|\xi|^{-2}t}|\widehat{u}_0|\notag\\
&\lesssim\chi_{\extt}(\xi)|\xi|^{-1}\mathrm{e}^{-c|\xi|^{-2}t}|\widehat{u}_0|,
\end{align}
where we used the mean value theorem to estimate the difference between two cosine functions, and the refined relations $\lambda_{\mathrm{I}}+|\xi|=O(|\xi|^{-3})$, $\lambda_{\mathrm{R}}+\frac{1}{2}|\xi|^{-2}=O(|\xi|^{-3})$ for $\xi\in\ml{Z}_{\extt}(N_0)$. Therefore, a combination of \eqref{Error-03} and \eqref{Error-04} leads to
\begin{align*}
\chi_{\extt}(\xi)\left|\widehat{u}-\cos(|\xi|t)\mathrm{e}^{-\frac{1}{2}|\xi|^{-2}t}\widehat{u}_0\right|\lesssim\chi_{\extt}(\xi)|\xi|^{-1}\mathrm{e}^{-c|\xi|^{-2}t}(|\widehat{u}_0|+|\widehat{u}_1|+|\xi|^{-2}|\widehat{v}_0|).
\end{align*}

\subsection{Decay estimate and asymptotic profile: Proof of Theorem \ref{Thm-Solution-Itself}}\label{Subsection-solution}
\hspace{5mm}To estimate the solution from the upper side, we apply the Plancherel theorem and the Hausdorff-Young inequality with the derived pointwise estimates \eqref{Est-u-INT}, \eqref{Est-u-BDD}, \eqref{Est-u-EXT} to obtain
\begin{align*}
\|u(t,\cdot)\|_{L^2}&\leqslant\|\chi_{\intt}(\xi)\widehat{u}(t,\xi)\|_{L^2}+\|\chi_{\bdd}(\xi)\widehat{u}(t,\xi)\|_{L^2}+\|\chi_{\extt}(\xi)\widehat{u}(t,\xi)\|_{L^2}\\
&\lesssim\left\|\chi_{\intt}(\xi)\mathrm{e}^{-c|\xi|^4t}\right\|_{L^2}\|(u_0,u_1,v_0)\|_{L^1\times L^1\times L^1}+\mathrm{e}^{-ct}\|(u_0,u_1,v_0)\|_{H^{\ell}\times H^{\ell-1}\times H^{\ell-3}}\\
&\quad+\left\|\chi_{\extt}(\xi)|\xi|^{-\ell}\mathrm{e}^{-c|\xi|^{-2}t}\right\|_{L^{\infty}}\|(u_0,u_1,v_0)\|_{H^{\ell}\times H^{\ell-1}\times H^{\ell-3}}\\
&\lesssim (1+t)^{-\frac{n}{8}}\|(u_0,u_1,v_0)\|_{L^1\times L^1\times L^1}+(1+t)^{-\frac{\ell}{2}}\|(u_0,u_1,v_0)\|_{H^{\ell}\times H^{\ell-1}\times H^{\ell-3}}
\end{align*}
with $\ell\geqslant0$, where we used \eqref{Technique-small} and \eqref{Technique-large} in the last line.

Let us recall the differential operators $\ml{G}^{(s)}(t,|D|)$ and $\ml{G}^{(l)}(t,|D|)$, whose symbols (or Fourier transforms) are given by
\begin{align*}
\widehat{\ml{G}}^{(s)}(t,|\xi|)=\chi_{\intt}(\xi)\mathrm{e}^{-|\xi|^4t}\ \ \mbox{and}\ \ \widehat{\ml{G}}^{(l)}(t,|\xi|)=\chi_{\extt}(\xi)\cos(|\xi|t)\mathrm{e}^{-\frac{1}{2}|\xi|^{-2}t}.
\end{align*}
From some derived estimates in the previous subsection, we may get the error estimate
\begin{align*}
&\left\|u(t,\cdot)-\ml{G}^{(s)}(t,|D|)\big(u_0(\cdot)-v_0(\cdot)\big)-\ml{G}^{(l)}(t,|D|)u_0(\cdot)\right\|_{L^2}\\
&=\left\|\chi_{\intt}(\xi)\big[\widehat{u}(t,\xi)-\mathrm{e}^{-|\xi|^4t}\big(\widehat{u}_0(\xi)-\widehat{v}_0(\xi)\big)\big]\right\|_{L^2}+\|\chi_{\bdd}(\xi)\widehat{u}(t,\xi)\|_{L^2}\\
&\quad+\left\|\chi_{\extt}(\xi)\big[\widehat{u}(t,\xi)-\cos(|\xi|t)\mathrm{e}^{-\frac{1}{2}|\xi|^{-2}t}\widehat{u}_0(\xi)\big]\right\|_{L^2}\\
&\lesssim\left\|\chi_{\intt}(\xi)\left(|\xi|^2\mathrm{e}^{-c|\xi|^4t}+\mathrm{e}^{-c|\xi|^2t}\right)\right\|_{L^2}\|(u_0,u_1,v_0)\|_{L^1\times L^1\times L^1}+\mathrm{e}^{-ct}\|(u_0,u_1,v_0)\|_{H^{\ell-1}\times H^{\ell-1}\times H^{\ell-3}}\\
&\quad+\left\|\chi_{\extt}(\xi)|\xi|^{-\ell}\mathrm{e}^{-c|\xi|^{-2}t}\right\|_{L^{\infty}}\|(u_0,u_1,v_0)\|_{H^{\ell-1}\times H^{\ell-1}\times H^{\ell-3}}\\
&\lesssim (1+t)^{-\min\left\{\frac{n}{8}+\frac{1}{2},\frac{n}{4}\right\}}\|(u_0,u_1,v_0)\|_{L^1\times L^1\times L^1}+(1+t)^{-\frac{\ell}{2}}\|(u_0,u_1,v_0)\|_{H^{\ell-1}\times H^{\ell-1}\times H^{\ell-3}}.
\end{align*}

Lastly, our next purpose is to derive the large time optimal estimate for the suitable higher regularity $\ell>\frac{n}{4}$ with $n\geqslant 1$. Clearly, \eqref{Est-05} implies the desired upper bound estimate in \eqref{Est-06} for large time. By the same manner as the one in \eqref{Est-Error-01}, we are able to get
\begin{align}
&\left\|\ml{G}^{(s)}(t,|D|)\big(u_0(\cdot)-v_0(\cdot)\big)-\ml{F}^{-1}_{\xi\to x}\left(\chi_{\intt}(\xi)\mathrm{e}^{-|\xi|^4t}\right)P_{u_0-v_0}\right\|_{L^2}\notag\\
&\lesssim t^{\frac{1}{8}}\left\|\chi_{\intt}(\xi)|\xi|\mathrm{e}^{-|\xi|^4t}\right\|_{L^2}\|(u_0,v_0)\|_{L^1\times L^1}+\left\|\chi_{\intt}(\xi)\mathrm{e}^{-|\xi|^4t}\right\|_{L^2}\int_{|y|\geqslant t^{\frac{1}{8}}}|u_0(y)-v_0(y)|\mathrm{d}y\notag\\
&\lesssim t^{\frac{1}{8}-\frac{n+2}{8}}\|(u_0,v_0)\|_{L^1\times L^1}+o(t^{-\frac{n}{8}})\notag\\
&=o(t^{-\frac{n}{8}})\label{Est-N-01}
\end{align}
for large time $t\gg1$, where we employed the $L^1$ integrability for $u_0-v_0$. Indeed, from the error estimate, we already derived
\begin{align*}
\left\|u(t,\cdot)-\ml{G}^{(s)}(t,|D|)\big(u_0(\cdot)-v_0(\cdot)\big)\right\|_{L^2}\lesssim t^{-\min\left\{\frac{n}{8}+\frac{1}{2},\frac{n}{4},\frac{\ell}{2} \right\}}\|(u_0,u_1,v_0)\|_{(H^{\ell}\cap L^1)\times (H^{\ell-1}\cap L^1)\times (H^{\ell-3}\cap L^1)}
\end{align*}
for large time $t\gg1$. Thanks to the condition $\ell>\frac{n}{4}$ so that $\min\left\{\frac{n}{8}+\frac{1}{2},\frac{n}{4},\frac{\ell}{2}\right\}>\frac{n}{8}$, we apply the triangle inequality to arrive at
\begin{align*}
\left\|u(t,\cdot)-\ml{F}^{-1}_{\xi\to x}\left(\chi_{\intt}(\xi)\mathrm{e}^{-|\xi|^4t}\right)P_{u_0-v_0}\right\|_{L^2}=o(t^{-\frac{n}{8}})
\end{align*}
and, analogously to \eqref{Est-07},
\begin{align*}
\|u(t,\cdot)\|_{L^2}\gtrsim t^{-\frac{n}{8}}|P_{u_0-v_0}|
\end{align*}
for large time $t\gg1$. Our proof is complete.

\appendix
\section{Large time behavior for the wave equation with the Riesz potential dissipation in $\mathbb{R}^n$}\label{Appendix-A}
\hspace{5mm}In this appendix, we briefly state the optimal estimate and the asymptotic profile as $t\gg1$ for the weakly dissipative wave equation \eqref{Wave-Riesz-Dissipation} in the whole space $\mb{R}^n$. Note that the well-posedness and some decay estimates for the Cauchy problem \eqref{Wave-Riesz-Dissipation} were studied by \cite[Theorems 2.1-2.3 with $\sigma=\delta=1$]{Said=2024}, where the restriction on the dimension $n$ can be removed due to the general definition of the Riesz potential instead of its explicit expression. Since the proof is similar to the one in Subsection \ref{Subsection-solution}, we omit some details. 
\begin{prop}\label{Prop-Riesz-wave}
	Let $\ell\geqslant0$. Suppose that the Cauchy data $w_0\in H^{\ell}\cap L^1$ and $w_1\in H^{\ell-1}\cap L^1$  for the weakly dissipative wave equation \eqref{Wave-Riesz-Dissipation}. Then, the solution satisfies the following decay estimate of the regularity-loss type:
	\begin{align*}
		\|w(t,\cdot)\|_{L^2}\lesssim (1+t)^{-\frac{n}{8}}\|(w_0,w_1)\|_{L^1\times L^1}+(1+t)^{-\frac{\ell}{2}}\|(w_0,w_1)\|_{H^{\ell}\times H^{\ell-1}}.
	\end{align*}
	Additionally, by assuming the suitable higher Sobolev regularity $\ell>\frac{n}{4}$ and $|P_{w_0}|\neq0$, the solution satisfies the following optimal estimate:
	\begin{align*}
		t^{-\frac{n}{8}}|P_{w_0}|\lesssim\|w(t,\cdot)\|_{L^2}\lesssim t^{-\frac{n}{8}}\|(w_0,w_1)\|_{(H^{\ell}\cap L^1)\times (H^{\ell-1}\cap L^1)}
	\end{align*}
	for large time $t\gg1$. Furthermore, considering the lower regular  Cauchy data $w_0\in H^{\ell-3}\cap L^1$, the solution satisfies the following refined estimate:
	\begin{align*}
		&\left\|w(t,\cdot)-\ml{G}^{(s)}(t,|D|)w_0(\cdot)-\ml{G}^{(l)}(t,|D|)w_0(\cdot)\right\|_{L^2}\notag\\
		&\lesssim (1+t)^{-\frac{n}{8}-\frac{1}{2}}\|(w_0,w_1)\|_{L^1\times L^1}+(1+t)^{-\frac{\ell}{2}}\|(w_0,w_1)\|_{H^{\ell-3}\times H^{\ell-1}}.
	\end{align*}
\end{prop}
\begin{remark}
The estimate of $w(t,\cdot)$ in the $L^2$ norm coincides with the one in \cite[Theorem 2.2 with $\sigma=\delta=1$]{Said=2024}, where we have removed the restrictions on the dimension $n$ by the general definition of the Riesz potential and on the additional regularity $\ell$ by the Bessel potential space (even with negative order). Moreover, we describe the large time profile of $w(t,\cdot)$ according to 
\begin{align*}
\left[\ml{G}^{(s)}(t,|D|)+\ml{G}^{(l)}(t,|D|) \right]w_0(\cdot)
\end{align*}
in the $L^2$ norm, as well as derive the optimal estimate when $\ell>\frac{n}{4}$. We believe these philosophies can be applied in the wave equation with general Riesz potential dissipations.
\end{remark}
\begin{proof}
By using the partial Fourier transform to the linear Cauchy problem \eqref{Wave-Riesz-Dissipation}, the solution $\widehat{w}=\widehat{w}(t,\xi)$ with the initial value $\widehat{w}_k=\widehat{w}_k(\xi)$ for $k=0,1$  is expressed by
\begin{align*}
\widehat{w}=\widehat{M}_0\widehat{w}_0+\widehat{M}_1\widehat{w}_1
\end{align*}
with the kernels in the Fourier space
\begin{align*}
\widehat{M}_0:=\frac{\lambda_1^{(w)}\mathrm{e}^{\lambda_2^{(w)}t}-\lambda_2^{(w)}\mathrm{e}^{\lambda_1^{(w)}t}}{\lambda_1^{(w)}-\lambda_2^{(w)}}\ \ \mbox{and}\ \ \widehat{M}_1:=\frac{\mathrm{e}^{\lambda_1^{(w)}t}-\mathrm{e}^{\lambda_2^{(w)}t}}{\lambda_1^{(w)}-\lambda_2^{(w)}},
\end{align*}
where the characteristic roots are formally given by
\begin{align*}
\lambda_{1,2}^{(w)}=\frac{-|\xi|^{-2}\pm\sqrt{|\xi|^{-4}-4|\xi|^2}}{2}.
\end{align*}
By using MacLaurin's formula $\sqrt{1+y}=1+\frac{y}{2}+O(y^2)$ for $|y|\ll1$ and the continuity of $\lambda_{1,2}^{(w)}$ with negative real parts for $\xi\in\ml{Z}_{\intt}(\varepsilon_0)\cup\ml{Z}_{\extt}(N_0)$ we know
\begin{itemize}
	\item for small frequencies: $\lambda_1^{(w)}=-|\xi|^4+O(|\xi|^{10})$ and $\lambda_2^{(w)}=-|\xi|^{-2}+O(|\xi|^4)$;
	\item for large frequencies: $\lambda_{1,2}^{(w)}=\lambda_{\mathrm{R}}^{(w)}\pm i\lambda_{\mathrm{I}}^{(w)}$ with $\lambda_{\mathrm{R}}^{(w)}=-\frac{1}{2}|\xi|^{-2}$ and $\lambda_{\mathrm{I}}^{(w)}=|\xi|+O(|\xi|^{-5})$;
	\item for bounded frequencies: $\mathrm{Re}\,\lambda_{1,2}^{(w)}<0$.
\end{itemize}
With the aid of these asymptotic expansions, we address
\begin{align*}
\chi_{\intt}(\xi)\widehat{M}_0&=\chi_{\intt}(\xi)\frac{[-|\xi|^4+O(|\xi|^{10})]\mathrm{e}^{-|\xi|^{-2}t+O(|\xi|^4)t}-[-|\xi|^{-2}+O(|\xi|^4)]\mathrm{e}^{-|\xi|^4t+O(|\xi|^{10})t}}{|\xi|^{-2}+O(|\xi|^4)},\\
\chi_{\intt}(\xi)\widehat{M}_1&=\chi_{\intt}(\xi)\frac{\mathrm{e}^{-|\xi|^4t+O(|\xi|^{10})t}-\mathrm{e}^{-|\xi|^{-2}t+O(|\xi|^4)t}}{|\xi|^{-2}+O(|\xi|^4)}.
\end{align*}
Applying the same approaches as those in Section \ref{Section-Solution}, we arrive at
\begin{align*}
\chi_{\intt}(\xi)|\widehat{w}|&\lesssim\chi_{\intt}(\xi)\mathrm{e}^{-c|\xi|^4t}(|\widehat{w}_0|+|\widehat{w}_1|),\\
\chi_{\intt}(\xi)|\widehat{w}-\mathrm{e}^{-|\xi|^4t}\widehat{w}_0|&\lesssim\chi_{\intt}(\xi) |\xi|^2\mathrm{e}^{-c|\xi|^4t}(|\widehat{w}_0|+|\widehat{w}_1|),
\end{align*}
where we noticed that $\mathrm{e}^{-|\xi|^{-2}t}\lesssim \mathrm{e}^{-ct}$ for $\xi\in\ml{Z}_{\intt}(\varepsilon_0)$. For large frequencies, thanks to the conjugate roots, our kernels can be rewrote by
\begin{align*}
\chi_{\extt}(\xi)\widehat{M}_0&=\chi_{\extt}(\xi)\left(\cos(\lambda_{\mathrm{I}}^{(w)}t)-\frac{\lambda_{\mathrm{R}}^{(w)}}{\lambda_{\mathrm{I}}^{(w)}}\sin(\lambda_{\mathrm{I}}^{(w)}t)\right)\mathrm{e}^{\lambda_{\mathrm{R}}^{(w)}t}\\
&=\chi_{\extt}(\xi)\left(\cos\left(|\xi|t+O(|\xi|^{-5})t\right)-\frac{-\frac{1}{2}|\xi|^{-2}}{|\xi|+O(|\xi|^{-5})}\sin\left(|\xi|t+O(|\xi|^{-5})t\right)\right)\mathrm{e}^{-\frac{1}{2}|\xi|^{-2}t}
\end{align*}
and
\begin{align*}
\chi_{\extt}(\xi)\widehat{M}_1&=\chi_{\extt}(\xi)\frac{\sin(\lambda_{\mathrm{I}}^{(w)}t)}{\lambda_{\mathrm{I}}^{(w)}}\mathrm{e}^{\lambda_{\mathrm{R}}^{(w)}t}\\
&=\chi_{\extt}(\xi)\frac{1}{|\xi|+O(|\xi|^{-5})}\sin\left(|\xi|t+O(|\xi|^{-5})t\right)\mathrm{e}^{-\frac{1}{2}|\xi|^{-2}t}.
\end{align*}
Thus, we may derive
\begin{align*}
	\chi_{\extt}(\xi)|\widehat{w}|&\lesssim\chi_{\extt}(\xi)\mathrm{e}^{-c|\xi|^{-2}t}(|\widehat{w}_0|+|\xi|^{-1}|\widehat{w}_1|),\\
	\chi_{\extt}(\xi)|\widehat{w}-\cos(|\xi|t)\mathrm{e}^{-\frac{1}{2}|\xi|^{-2}t}\widehat{w}_0|&\lesssim\chi_{\extt}(\xi) \mathrm{e}^{-c|\xi|^{-2}t}(|\xi|^{-3}|\widehat{w}_0|+|\xi|^{-1}|\widehat{w}_1|).
\end{align*}
Let us recall the tools \eqref{Technique-small} as well as \eqref{Technique-large}. Due to an exponential decay estimate for bounded frequencies, making use of the derived pointwise estimates in the Fourier space, we obtain
\begin{align*}
\|w(t,\cdot)\|_{L^2}&\lesssim\left\|\chi_{\intt}(\xi)\mathrm{e}^{-c|\xi|^4t}\right\|_{L^2}\|(w_0,w_1)\|_{L^1\times L^1}+\mathrm{e}^{-ct}\|(w_0,w_1)\|_{H^{\ell}\times H^{\ell-1}}\\
&\quad+\sup\limits_{|\xi|\geqslant N_0}\left(|\xi|^{-\ell}\mathrm{e}^{-c|\xi|^{-2}t}\right)\|(w_0,w_1)\|_{H^{\ell}\times H^{\ell-1}},
\end{align*}
which implies the desired estimate for the solution $w$, and
\begin{align*}
&\left\|w(t,\cdot)-\ml{G}^{(s)}(t,|D|)w_0(\cdot)-\ml{G}^{(l)}(t,|D|)w_0(\cdot)\right\|_{L^2}\\
&\lesssim\left\|\chi_{\intt}(\xi)|\xi|^2\mathrm{e}^{-c|\xi|^4t}\right\|_{L^2}\|(w_0,w_1)\|_{L^1\times L^1}+\mathrm{e}^{-ct}\|(w_0,w_1)\|_{H^{\ell-3}\times H^{\ell-1}}\\
&\quad+\sup\limits_{|\xi|\geqslant N_0}\left(|\xi|^{-\ell}\mathrm{e}^{-c|\xi|^{-2}t}\right)\|(w_0,w_1)\|_{H^{\ell-3}\times H^{\ell-1}},
\end{align*}
which shows the desired error estimate. Finally, for the higher regularity $\ell>\frac{n}{4}$ with $n\geqslant 1$, replacing $u_0-v_0$ in \eqref{Est-N-01} by $w_0$, and using 
\begin{align*}
\|w(t,\cdot)-\ml{G}^{(s)}(t,|D|)w_0(\cdot)\|_{L^2}\lesssim t^{-\min\left\{\frac{n}{8}+\frac{1}{2},\frac{\ell}{2}\right\}}\|(w_0,w_1)\|_{(H^{\ell}\cap L^1)\times (H^{\ell-1}\cap L^1)}
\end{align*}
for large time $t\gg1$, we complete the desired optimal lower bound estimate by the same way as the one for $u$ in Subsection \ref{Subsection-solution}. The proof is finished.
\end{proof}

\section*{Acknowledgments} 
 Wenhui Chen is supported by the National Natural Science Foundation of China (grant No. 12301270, grant No. 12171317), 2024 Basic and Applied Basic Research Topic--Young Doctor Set Sail Project (grant No. 2024A04J0016). Yan Liu is supported by Guangdong
 Basic and Applied Basic Research Foundation (Grant 2023A1515012044).

% ------------------------------------------------------------------------
\end{document}